\documentclass[smallextended]{svjour3}       % onecolumn (second format)
\smartqed  % flush right qed marks, e.g. at end of proof
\usepackage{graphicx}
\usepackage{amsmath}
\DeclareMathOperator{\sgn}{sgn}
 
\begin{document}

\title{Semi-implicit methods for advection equations with explicit forms of numerical solution. 
\thanks{The work was supported by the grant VEGA 1/0709/19 and APVV-19-0460.}
}
%\subtitle{The title is very preliminary ...}

\titlerunning{Semi-implicit methods for advection equations} 

\author{Peter Frolkovi\v{c}         \and Svetlana Kri\v{s}kov\'a 
\and Michaela Rohov\'a \and Michal \v{Z}erav\'y
}


\institute{Peter Frolkovi\v{c} \and Svetlana Kri\v{s}kov\'a  \and Michaela Rohov\'a \and   Michal \v{Z}erav\'y \at
              Dept. of Mathematics and Descriptive Geometry, Faculty of Engineering STU\\
              Radlinsk\'eho 11, 810 05 Bratislava, Slovakia \\
              \email{peter.frolkovic@stuba.com}           %  \\
}

\date{Received: \today / Accepted: The correct dates will be entered by the editor}

\maketitle

\begin{abstract}
We present a parametric family of semi-implicit second order accurate numerical methods for non-conservative and conservative advection equation for which the numerical solutions can be obtained in a fixed number of forward and backward alternating substitutions. The methods use a novel combination of implicit and explicit time discretizations for one-dimensional case and the Strang splitting method in several dimensional case. The methods are described for advection equations with a continuous variable velocity that can change its sign inside of computational domain. The methods are unconditionally stable in the non-conservative case for variable velocity and for variable numerical parameter. Several numerical experiments confirm the advantages of presented methods including an involvement of differential programming to find optimized values of the variable numerical parameter. 

\end{abstract}

\section{Introduction}
\label{intro}

Implicit (or semi-implicit) numerical schemes are useful numerical methods to solve advection dominated problems in several circumstances \cite{arbogast2020third,carciopolo2019conservative,frolkovivc2016numerical,frolkovivc2018semi,fuhrmann2001stability,hadjimichael2021positivity,hahn2019iterative,ibolya2020numerical,knodel2020global,mikula2014inflow,partl2016numerical,polivka2014compositional,puppo2021quinpi,qin2018implicit,radu2010analysis,zhang2021numerical}. They can avoid or reduce significantly the main disadvantage of fully explicit schemes that are implemented on a fixed mesh with a finite stencil in numerical discretization. Such explicit schemes based on finite difference or finite volume methods require a stability restriction on the choice of discretization parameters \cite{leveque2002finite}, e.g. the time step. This can be disadvantageous in several situations like highly variable velocity field or nonuniform space discretization step \cite{frolkovivc2018semi,may2017explicit} or long time simulations when approaching stationary solutions \cite{qin2018implicit}. Another important type of problems when implicit schemes can be useful are stiff differential equations involving the advection term \cite{frolkovivc2016numerical,izzo2017highly,knodel2020global,partl2016numerical}.

In an ideal case, the implicit and semi-implicit methods can offer an unconditional stability that make them convenient tool to solve numerically the problems having previously mentioned complications. The price to pay is that the numerical solution must be obtained in general by solving linear algebraic system of equations. Our aim is to propose a semi-implicit second order accurate numerical method having unconditional stability where the obtained linear systems can be solved in a small given number of alternating substitutions. Consequently, such scheme can compete well with fully explicit schemes that themselves lack the advantage of unconditional stability.

In \cite{frolkovivc2018semi} a novel parametric family of semi-implicit second order accurate numerical schemes for linear non-conservative advection is introduced. We base our work here on this result that is extended in several aspects. We treat here only representative advection equations, but the derived numerical schemes can be used for more involved advection dominated problems. 

Additionally to \cite{frolkovivc2018semi}, we show how one-dimensional discretized problems with a velocity that can change its sign several times in domain can be solved by one forward and one backward substitution as known for so called  fast sweeping methods \cite{lozano2021implicit,zhao2005fast}. The derivation of second order scheme is based here on the discretization of an error term for the fully implicit first order accurate upwind scheme. Consequently, it is easy to apply some limiter procedures based on a blending between these two types of schemes if necessary. The free parameter of this family of schemes is introduced in a more convenient way than in \cite{frolkovivc2018semi}. The parameter can vary in space and time that we demonstrate in one numerical experiment by letting find its optimized values using differential programming tools.

In this work we present for the first time the semi-implicit scheme also for the linear advection equation in the conservative form. The scheme is based on finite volume method and it offers an exact local mass balance property in a discrete form that is confirmed also by numerical experiments. Nevertheless, to insure such property, the local variability of the parameter in the method is not allowed and the unconditional von Neumann stability can be shown only for constant velocity case. Therefore, if the advection dominated problem can be expressed in a non-conservative form, e.g. for the divergence free velocity, we prefer the more flexible semi-implicit non-conservative method. As confirmed by several numerical experiments, the errors in mass balance for chosen numerical examples diminish fast with a mesh refinement. 

Finally, we apply the derived one-dimensional non-conservative scheme for two-dimensional advection equation using Strang splitting \cite{leveque2002finite,uccar2019operator} that preserves the second order accuracy. Such treatment requires to solve only a given number of one-dimensional problems in alternating directions that can be solved by our proposed numerical scheme. One can show that such splitting scheme preserves the property of unconditional stability \cite{uccar2019operator} that is not the case of unsplitted version presented in \cite{frolkovivc2018semi}. 

Although we describe the semi-implicit schemes only for simple squared domain with uniform mesh, we note that its extension to complex three-dimensional domains with unstructured polyhedral mesh is published, e.g., in \cite{hahn2019iterative}. Moreover, the uniform grids can be used for nontrivial two-dimensional domains with unfitted grids as described and illustrated in \cite{frolkovivc2018semi} where the unconditional stability of the scheme is very important due to the presence of arbitrary small size of so-called cut cells \cite{may2017explicit}. Furthermore, we treat here only problems with smooth solutions that is relevant for level set methods, but we can quote first preliminary results of extensions for this type of schemes also for nonlinear conservation laws \cite{zeravy2021} with discontinuous solutions that is content of our future research.

We proceed as follows. In Section \ref{sec-nonc} we derive the parametric family of one-dimensional semi-implicit non-conservative second order accurate advection scheme. The properties of the scheme are formulated and proved in two propositions. In Section \ref{sec-c} we derive analogous conservative advection schemes. In Section \ref{sec-2d} we introduce the Strang splitting for problems in two dimension that preserves the accuracy and the unconditional stability of 1D scheme. In Section \ref{sec-ne} we illustrate all properties on several numerical examples. 

\section{Non-conservative advection equation}
\label{sec-nonc}

The linear advection equation in the non-conservative form is given by
\begin{equation}
    \label{1dadv}
    \partial_t \phi + v \partial_x \phi = 0 \,, \quad \phi(x,0) = \phi^0(x) \,,
\end{equation}
where the velocity function $v=v(x)$ is a given continuous function. The unknown function $\phi=\phi(x,t)$ is prescribed at $t=0$ by the given function $\phi^0=\phi^0(x)$ and $\phi(x,t)$ should be determined for $t>0$ and $x \in (0,L) \subset R$. The boundary values are prescribed by given functions $\phi_0=\phi_0(t)$ and $\phi_L=\phi_L(t)$ only if an inflow situation occurs at the boundary, namely,
\begin{equation}
\label{1dbc}
\phi(0,t) = \phi_0(t) \,, \,\, \hbox{if   } v(0) \ge 0 \,, \quad
\phi(L,t) = \phi_L(t) \,, \,\,  \hbox{if  } v(L) \le 0 \,.
\end{equation}

The model equation (\ref{1dadv}) can be solved by finding so called characteristic curves, the characteristics, generated by $v$ and using the fact that the solution $\phi$ is constant along characteristics \cite{leveque2002finite}.  

In what follows, we use the following common notations of finite difference methods. We denote $x_i = i h$, $i=0,1,\ldots,I$ for a chosen $I$, where $h=L/I$, and $t^n=n \tau$, $n=0,1,\ldots,N$ for a chosen $N$ and $\tau>0$. Our aim is to find the approximations $\phi_i^n$ of $\phi(x_i,t^n)$. The initial values are given by $\phi_i^0 = \phi^0(x_i)$, and for $i=0$ and $i=I$ they can be replaced by the values using the boundary conditions (\ref{1dbc}).

\begin{remark}
\label{rem1}
To simplify our presentation, we consider a particular form of the velocity $v$ being a piecewise linear function of $x$ that can be determined by its discrete values $v_i = v(x_i)$ for each mesh. Such situation is typical if the velocity field is obtained by a numerical solution of some flow equation or when the velocity function is approximated from its values in mesh points. We allow that the velocity can change its sign inside of $(0,L)$ when the points for which $v=0$ can be then easily determined from its piecewise linear form. 

Later we distinguish the case of so called expanding characteristics for points $\bar x \in (x_i,x_{i+1})$ such that $v(\bar x) =0$ and $v'(\bar x) > 0$. In this case, the regions $x \le \bar x$ and $x \ge \bar x$  are "separated" from each other that may require a special treatment in numerical schemes, see later. Note that if the velocity $v$ in (\ref{1dadv}) depends also on $t$, we use an approximation by fixing its value on each time subinterval $[t^n,t^{n+1}]$ using a representative value, e.g., at $t=t^n+\tau/2$. In such case, the positions of zero points $\bar x$ have to be redetermined in each time subinterval.
\end{remark}

\subsection{First order accurate fully implicit scheme}
\label{sec-1st}

The scheme can be derived using the backward Euler method for the time discretization and the upwind one-sided finite difference for the space discretization,
\begin{equation}
%\label{1d1o}
\nonumber
\frac{\phi_i^{n+1}-\phi_i^n}{\tau} + v_i^+ \, \frac{\phi_i^{n+1}-\phi_{i- 1}^{n+1}}{h} + v_i^- \, \frac{\phi_{i+1}^{n+1}-\phi_{i}^{n+1}}{h} = 0 ,
\end{equation}
where $v^+ := \max\{0,v\}$ and $v^- := \min\{0,v\}$.  
Denoting (signed) Courant numbers at mesh points
\begin{equation}
%\label{cn}
C_i := \frac{\tau v_i}{h} \,,
\end{equation}
we can write the scheme in the concise form
\begin{equation}
\label{1d1or}
\phi_i^{n+1} = \left( 1 + |C_i|\right)^{-1} \left(\phi_i^n +  |C_i|\phi_{i\mp 1}^{n+1} \right)
\end{equation}
with $\mp = -\sgn(C_i)$. The scheme (\ref{1d1or}) can be used for $n=0,1,\ldots,N$ and for the inner nodes with $i=1,2,\ldots,I-1$. For boundary nodes the scheme can be applied only if $C_0<0$ or $C_I>0$, otherwise the Dirichlet boundary conditions (\ref{1dbc}) must be used,
\begin{eqnarray}
\label{1dbcn}
\phi_0^{n+1} = \phi_0(t^{n+1}) \,, \,\, \hbox{if   } C_0 \ge 0 \,, \\[1ex]
\label{1dbcnB}
\phi_I^{n+1} = \phi_L(t^{n+1}) \,, \,\,  \hbox{if  } C_I \le 0 \,.
\end{eqnarray}
% or one simplifies (\ref{1d1or}) for the outflow boundaries to
% \begin{eqnarray}
% \label{bc0}
% \phi_0^{n+1} = \left( 1 - C_0 \right)^{-1}\left(\phi_0^n  - C_0\phi_{1}^{n+1}\right) \,, \,\,
% \hbox{if  } C_0 < 0 \,, \\[1ex]
% \label{bcL}
% \phi_I^{n+1} = \left( 1 + C_I \right)^{-1}\left(\phi_I^n +  C_I\phi_{I- 1}^{n+1}\right) \,, \,\,
% \hbox{if  } C_I > 0 \,.
% \end{eqnarray}

Finally, we modify (\ref{1d1or}) for the special case of expanding characteristics as described in Remark \ref{rem1}. Namely, if there exists $i^*$ such that
\begin{eqnarray}
\label{sign}
    v_{i^*} < 0 \hbox{  and  } v_{1+i^*} > 0 \,,
\end{eqnarray}
then there is a point $\bar x_{i^*} \in (x_{i^*}, x_{1+i^*})$ where the piecewise linear velocity attains the zero value.
For such point with expanding characteristics we decouple the computations of $\phi_{i^*}^{n+1}$ and $\phi_{1+i^*}^{n+1}$ from each other.
Namely, the scheme (\ref{1d1or}) is replaced for $i=i^*$ and $i=1+i^*$ by the explicit definitions
\begin{eqnarray}
\label{zero}
\begin{array}{l}
\phi_{i^*}^{n+1} = \left( 1  - C_{i^*}\right)^{-1}\left(\phi_{i^*}^n  - C_{i^*} \phi_{i^{**}}^{n}\right)\,, \\[1ex]    
    \phi_{1+i^*}^{n+1} = \left( 1 + C_{1+i^*}\right)^{-1}\left(\phi_{1+i^*}^n +  C_{1+i^*} \phi_{i^{**}}^{n}\right) \,.
\end{array}
\end{eqnarray}
To use it, we determine the location of $\bar x_{i^*}$ from the linear interpolation, 
\begin{equation}
%    \label{li}
\nonumber
    \bar x_{i^*} = x_{i^*} + \frac{v_{i^*}}{v_{i^*} - v_{1+i^*}} h \,.
\end{equation}
Afterwards, we define the interpolated value $\phi_{i^{**}}^n \approx \phi(\bar x_{i^*},t^n)$ by
\begin{equation}
%    \label{phiplus}
\nonumber
    \phi_{i^{**}}^n := \phi_{i^*}^n + \frac{\bar x_{i^*}-x_{i^*}}{h} (\phi_{1+i^*}^n - \phi_{i^*}^n)
\end{equation}
and 
\begin{equation}
%    \label{cnspecial}
\nonumber
C_{1+i^*} := \frac{\tau v_{1+i^*}}{x_{1+i^*}-\bar x_{i^*}} \,, \quad
C_{i^*} := 
  \frac{\tau v_{i^*}}{\bar x_{i^*}-x_i} \,.
\end{equation}
The schemes (\ref{zero}) are then determined using the approximation that $\phi(x_{i^*},t) \equiv \phi_{i^{**}}^n$ for $t \in [t^n,t^{n+1}]$ as $v(\bar x_{i^*})=0$.

Formally, the schemes (\ref{1d1or}) - (\ref{zero}) represent a system of linear algebraic equations. Nevertheless, each $i$-th equation in (\ref{1d1or}) - (\ref{zero}) contains at most one neighbor value, either $\phi_{i-1}^{n+1}$ or $\phi_{i+1}^{n+1}$ (or none). Such linear equations can be solved
using one forward and backward substitution as known for fast sweeping methods \cite{zhao2005fast} or, alternatively, using fractional time steps as described in \cite{lozano2021implicit}. The latter method will be applied and described in the next section.

% Alternatively, one can formulate the substitutions analogously to \cite{lozano2021implicit} when the discrete equations are split into two fractional steps that must be applied in the defined order for $i$, namely
% \begin{eqnarray}
% \label{1d1orf1}
% \phi_i^{n+1/2} = \left( 1 + C_i^+\right)^{-1} \left(\phi_i^n +  C_i^+\phi_{i- 1}^{n+1/2} \right) \,, 
% \,\, i=0,1,\ldots,I\\[1ex]
% \label{1d1orf2}
% \phi_i^{n+1} = \left( 1 - C_i^-\right)^{-1} \left(\phi_i^{n+1/2} - C_i^-\phi_{i+1}^{n+1}\right) \,,
% \,\, i=I,I-1,\ldots,0 \,.
% \end{eqnarray}
% The boundary conditions have to be appropriately included, namely taking (\ref{1dbcn}) instead of (\ref{1d1orf1}) for $i=0$ if $C_0 \ge 0$, and (\ref{1dbcnB}) instead of (\ref{1d1orf2}) for $i=I$ if $C_I \le 0$.

We can now summarize the advantages of the first order accurate implicit upwind scheme. The numerical solution is defined by (\ref{1d1or}) and (\ref{zero}) as convex combinations of the neighbor values for arbitrary large Courant numbers. Therefore, the scheme is unconditionally stable and it insures a discrete minimum and maximum principle for any choice of $h$ and $\tau$. This is especially critical if (\ref{zero}) has to be used as the corresponding Courant numbers $C_{i^*}$ or $C_{1+i^*}$ can be arbitrary large.  

The main disadvantage is the low accuracy of the method that can be demonstrated for many test examples of practical interest. Therefore, we derive an extension of this method in a form of the second order accurate  semi-implicit methods for which the numerical solution can be again determined explicitly.

\subsection{Parametric family of second order accurate semi-implicit schemes}
\label{sec-2nd}

We begin by presenting a form of the one-step error \cite{leveque2002finite} of the first order accurate scheme (\ref{1d1or}). To do so, we express the values $\phi(x_i,t^n)$ and $\phi(x_{i\pm 1},t^{n+1})$ using Taylor series at $(x_i,t^{n+1})$. First, we have for $t=t^{n+1}$
\begin{eqnarray}
    \label{taylorx}
    \phi(x_{i \pm 1},\cdot) = \phi(x_i,\cdot) \pm h \partial_x \phi(x_i,\cdot) + \frac{h^2}{2} \partial_{xx} \phi(x_i,\cdot) + \mathcal{O}(h^3)  \,.
\end{eqnarray}
Second, we have for $x=x_i$
\begin{eqnarray}
    \label{taylor}
    \phi(\cdot,t^n) = \phi(\cdot,t^{n+1}) - \tau \partial_t \phi(\cdot,t^{n+1}) + \frac{\tau^2}{2} \partial_{tt} \phi(\cdot,t^{n+1}) + \mathcal{O}(\tau^3) = \quad \\[1ex] \nonumber
    \phi(\cdot,t^{n+1}) + \tau v_i \partial_x \phi(\cdot,t^{n+1}) - \frac{\tau^2}{2}  v_i \partial_{tx} \phi(\cdot,t^{n+1}) + \mathcal{O}(\tau^3) \,,
\end{eqnarray}
where one exploits that $\partial_t \phi=-v \partial_x \phi$ and $\partial_{tt} \phi= - v \partial_{tx} \phi$. Such approach is often called Lax-Wendroff procedure  \cite{leveque2002finite}. Opposite to its standard form when all time derivatives in (\ref{taylor}) are replaced by spatial derivatives using the equation (\ref{1dadv}),  we allow also the mixed derivatives in (\ref{taylor}).

Using (\ref{taylor}) and (\ref{taylorx}) we obtain
that the discrete values $\phi(x_i,t^{n+1})$ of exact solution fulfill the first order accurate scheme (\ref{1d1or}) with the leading term of the one-step error
\begin{equation}
    \label{error}
    e_2(x_i,t^{n+1})  := - \frac{\tau^2}{2} v_i \partial_{tx} \phi(x_i,t^{n+1}) + \frac{h \tau}{2} |v_i| \partial_{xx} \phi(x_i,t^{n+1})  \,.
\end{equation}

Now, to extend the scheme (\ref{1d1or}) to be second order accurate, we have to approximate the derivatives in the one-step error $e_2$ with at least first order accurate approximations. Our aim is to derive a parametric family of semi-implicit upwind schemes that has a convenient stencil in its implicit part and that is unconditionally stable \cite{frolkovivc2018semi}.

To do so we define the parametric ``upwind based'' approximations of $\partial_x \phi(x_i,t^*)$ for $*=n$ or $*=n+1$,
\begin{equation}
\label{kappa}
\begin{array}{l}
h \, \partial_x^{\alpha-} \phi_i^* := \alpha (\phi_{i}^* - \phi_{i-1}^*) + (1-\alpha) (\phi_{i+1}^* - \phi_{i}^*) \\[1ex]
h \, \partial_x^{\alpha+} \phi_i^* := \alpha (\phi_{i+1}^* - \phi_{i}^*) + (1-\alpha) (\phi_{i}^* - \phi_{i-1}^*) \,.
\end{array}
\end{equation}
The approximation $\partial_x^{\alpha-} \phi_i^*$ will be used in (\ref{error}) if $v_i>0$ and $\partial_x^{\alpha+} \phi_i^*$ if $v_i<0$. Additionally, applying the standard backward finite difference for the time derivative, and the upwind finite difference  for the spatial derivatives, we obtain
\begin{eqnarray}
\nonumber
\begin{array}{l}
e_2(x_i,t^{n+1}) \approx 
- \dfrac{\tau}{2} (v_i^+ \partial_x^{\alpha-} + v_i^- \partial_x^{\alpha+}) (\phi_i^{n+1} -  \phi_i^{n}) \quad \\[1ex]
+ \dfrac{\tau}{2} v_i^+ \partial_{x}^{\alpha-} \left(  \phi_i^{n+1} -  \phi_{i-1}^{n+1}\right)  -   \dfrac{\tau}{2} v_i^- \partial_{x}^{\alpha+} \left(  \phi_{i+1}^{n+1} - \phi_i^{n+1}\right) \,.
\end{array}
\end{eqnarray}
Doing it this way we see that the terms with $\partial_x^{\alpha \mp} \phi_i^{n+1}$ cancel. 

Combining now the first order and the second order accurate approximations, we obtain the final semi-implicit scheme of the form
\begin{eqnarray}
\label{1d2obetween}
\phi_i^{n+1} + C_i^+  \left( \phi_i^{n+1} - \phi_{i-1}^{n+1} - \frac{h}{2} \partial_x^{\alpha-} \phi_{i-1}^{n+1} \right) + \\[1ex]\nonumber
 C_i^-  \left( \phi_{i+1}^{n+1}   - \phi_{i}^{n+1} -  \frac{h}{2} \partial_x^{\alpha+} \phi_{i+1}^{n+1} \right) = \phi_i^n - C_i^+ \frac{h}{2} \partial_x^{\alpha-} \phi_{i}^{n}
 - C_i^- \frac{h}{2} \partial_x^{\alpha+} \phi_{i}^{n} \,,
\end{eqnarray}
where one can clearly distinguish between the contributions of two approximations having the different order of accuracy. One can write (\ref{1d2obetween}) in the concise form using $\mp = - \sgn(C_i)$,
\begin{eqnarray}
\label{1d2obetweenelegant}
\phi_i^{n+1} + |C_i|  \left( \phi_i^{n+1} - \phi_{i\mp 1}^{n+1} \mp \frac{h}{2} \partial_x^{\alpha \mp} \phi_{i \mp 1}^{n+1} \right) = \phi_i^n \mp  |C_i| \frac{h}{2} \partial_x^{\alpha \mp} \phi_{i}^{n} \,.
\end{eqnarray}

The value $\alpha=0$ choose ``downwind'' one-sided finite difference in (\ref{kappa}) and $\alpha=1$ the upwind one. The case $\alpha = 0.5$ results in the central finite difference in (\ref{kappa}). In general, the parameters $\alpha$ in (\ref{kappa}) can be different for each $i$ and $n$. 

We summarize the accuracy and stability properties of the scheme in the following Proposition. 
% Similar accuracy (or consistency) and stability requirements are used for numerical schemes to solve linear PDEs are sufficient ingredients for a convergence of related numerical methods according to Lax equivalence theorem \cite{leveque2002finite}.

\begin{proposition}
\label{propacc}
Let $\phi$ be a smooth solution of the linear advection equation  (\ref{1dadv}), then the one-step error of the scheme (\ref{1d2obetweenelegant}) is given by
\begin{eqnarray}
    \label{trunc}
    e_3(x_i,t^{n+1}) = \\
    \nonumber \frac{1}{12} C_i h \left( \tau^2 \partial_{ttx} \mp 3 (1- 2\alpha_i) \tau h \partial_{txx} + 2 (1-3 \alpha_i) h^2 \partial_{xxx} \right) \phi(x_i,t^{n+1}) \,.
\end{eqnarray}
Consequently, the scheme is $2^{nd}$ order accurate and  in the case of constant velocity $v$  it is $3^{rd}$ order accurate if
\begin{equation}
    \label{kappav}
   \alpha_i = \frac{2 + |C_i|}{6}\,. 
\end{equation}
The scheme is unconditionally stable in the sense of von Neumann stability analysis for any $\alpha_i \ge 0$. 

\end{proposition}

\begin{proof}
The one step error (\ref{trunc}) can be obtained analogously as described in this section for the derivation of the one-step error $e_2$ for the $1^{st}$ order scheme. 
To prove the $3^{rd}$ order accuracy in the case of constant velocity we use $C_i \equiv C = v h / \tau$ and $\partial_{txx} \phi = - v \partial_{xxx} \phi$ and $\partial_{ttx} \phi = - v \partial_{txx} \phi = v^2 \partial_{xxx} \phi$, when 
\begin{eqnarray}
    \label{truncp}
    e_3(x_i,t^{n+1}) = \\
    \nonumber \frac{1}{12} C h \left( \tau^2 v^2 \pm 3 (1 - 2\alpha_i) \tau h v + 2 (1-3 \alpha_i) h^2 \right) \partial_{xxx} \phi(x_i,t^{n+1}) = \\[1ex]
    \nonumber \frac{1}{12} C h^3 \left( C^2  + 3 |C| - 6 \alpha_i |C| + 2 - 6 \alpha_i) \right) \partial_{xxx} \phi(x_i,t^{n+1}) = \\[1ex]
    \nonumber \frac{1}{12} C h^3 (|C| + 1) \left(  |C| + 2 - 6 \alpha_i \right) \partial_{xxx} \phi(x_i,t^{n+1}) \,.
\end{eqnarray}
We remind that $\pm = \sgn(v)$ and $\mp = -\sgn(v)$. Consequently for the choice (\ref{kappav}) the error term $e_3$ in (\ref{truncp}) vanishes.

Next we prove the linear stability of (\ref{1d2obetweenelegant}) using von Neumann stability analysis \cite{leveque2002finite,frolkovivc2018semi,arbogast2020third}. To do so we consider $C_i>0$ (the other case is treated analogously) and we rewrite (\ref{1d2obetweenelegant}) with short notation $C=C_i$ and $\alpha=\alpha_i$ to
\begin{eqnarray}
\label{1d2obetweenelegant000}
\phi_i^{n+1} + C  \left( \phi_i^{n+1} - \phi_{i-1}^{n+1} - \frac{1}{2} \left((1-\alpha)\phi_{i}^{n+1}+(2\alpha-1)\phi_{i-1}^{n+1}-\alpha\phi_{i-2}^{n+1} \right) \right) \\[1ex]\nonumber
= \phi_i^n -  C \frac{1}{2} \left((1-\alpha)\phi_{i+1}^{n}+(2\alpha-1)\phi_{i}^{n}-\alpha\phi_{i-1}^{n}\right) \,.
\end{eqnarray}
Next, we consider discrete Fourier modes in the complex plane
\begin{equation}
\nonumber
    \epsilon_{i+ k}^n := \exp(-\lambda t^n) \exp(\imath k \theta) \,, \quad \theta \in [-\pi,\pi] \,, \,\, k=-2,-1,\ldots,1 \,.
\end{equation}
Let $S:=\exp(-\lambda \tau)$  and we use following straightforward relations
\begin{equation}
    \label{eps}
\epsilon_{i}^{n+1} = S \epsilon_{i}^n \,, \quad  \epsilon_{i+k}^n = \left(\cos(k
\theta) + \imath \sin(k \theta) \right) \epsilon_i^n \,.
\end{equation}
Our aim is to show that $|S|\le 1$. Replacing $u$ in (\ref{1d2obetweenelegant000}) by $\epsilon$ \cite{leveque2002finite,frolkovivc2018semi,arbogast2020third} and using (\ref{eps}) we obtain after simple algebraic manipulations
% \begin{eqnarray}
% \nonumber
% S \left(1 + C_i \left( 1 - E(-\theta) -
% \frac{1}{2} \left(1-\alpha+(2\alpha-1) E(-\theta)-\alpha E(- 2 \theta) \right) \right) \right) \\[1ex] 
% \nonumber =
% \left(1 - C_i \frac{1}{2} \left( (1-\alpha) E(\theta)+2\alpha-1-\alpha E(- \theta) \right)  \right)\end{eqnarray}
% Clearly, the factor $S$ must fulfill
% \begin{eqnarray}
% \label{S}
% S  =
% \frac{1 +  \frac{C_i}{2} \left( 1-2\alpha + (\alpha-1) E(\theta)+\alpha E(-\theta) \right)}{1 + \frac{C_i}{2}
% \left(1+\alpha-(2\alpha+1) E(-\theta)+\alpha E(-2 \theta) \right)} \end{eqnarray}
\begin{eqnarray}
\nonumber
S \left(1 + \frac{C}{2}\left( 1+\alpha-(2\alpha+1)(\cos(\theta)-\imath \sin(\theta))+\alpha (\cos(2 \theta)-\imath \sin(2 \theta)) \right) \right) = \\[1ex]
\nonumber
1 + \frac{C}{2} \left(1-2\alpha+(2\alpha-1)\cos(\theta) - \imath \sin(\theta) \right) \,.
\end{eqnarray}
Solving the last relation as a complex algebraic equation for the real and imaginary part of $S$ and computing $|S|^2$ from the result using Mathematica \cite{Mathematica} we obtain that $|S|^2 = \gamma/\delta$
with
\begin{eqnarray}
\nonumber
\gamma = (\alpha -1) \alpha  C^2 \cos (2 \theta )-(2
   \alpha -1) C \left((2 \alpha
   -1) C-2\right)\cos (\theta ) + \\[1ex]\nonumber
   C \left(\alpha  \left(3
   (\alpha -1) C-4\right)+C+2\right)+2\\[1ex]
   \nonumber
\delta = \alpha  C \left(\alpha 
   C+C+2\right) \cos (2 \theta ) - (2 \alpha +1) C \left(2 \alpha  C+C+2\right) \cos
   (\theta ) + \\[1ex] \nonumber
   C \left(\alpha  \left(3 (\alpha +1)
   C+2\right)+C+2\right)+2
% \gamma = (\alpha -1) \alpha  C_i^2 \cos (2 \theta )-(2
%   \alpha -1) C_i \left((2 \alpha
%   -1) C_i-2\right)\cos (\theta ) + \\[1ex]\nonumber
%   C_i \left(\alpha  \left(3
%   (\alpha -1) C_i-4\right)+C_i+2\right)+2\\[1ex]
%   \nonumber
% \delta = \alpha  C_i \left(\alpha 
%   C_i+C_i+2\right) \cos (2 \theta ) - (2 \alpha +1) C_i \left(2 \alpha  C_i+C_i+2\right) \cos
%   (\theta ) + \\[1ex] \nonumber
%   C_i \left(\alpha  \left(3 (\alpha +1)
%   C_i+2\right)+C_i+2\right)+2
\end{eqnarray}
From the definition one has $\gamma\ge 0$ and $\delta\ge 0$. Comparing the nominator $\gamma$ and denominator $\delta$ we get
\begin{eqnarray}
\nonumber
\delta - \gamma = \alpha C (1 + C) \left( 2 \cos(2 \theta) - 8 \cos(\theta) + 6\right)
= 16 \alpha C (1 + C) \sin^4(\theta/2)
\end{eqnarray}
Clearly, for $\alpha \ge 0$ and $C \ge 0$ we have that $|S|^2 \le 1$, therefore the scheme is unconditionally stable using the von Neumann stability analysis. Analogous results is obtained in the case $C_i \le 0$.
\end{proof}

\begin{proposition}
\label{propexpl}
The numerical solution of linear advection equation (\ref{1dadv}) using (\ref{1d2obetweenelegant}) can be obtained explicitly by one forward and one backward substitution:
\begin{eqnarray}
\label{1d2o_kappaelegant1}
\begin{array}{l}
\phi_i^{n+1/2} = \dfrac{2 \phi_i^n + C_i^+ \left((1+2 \alpha_i) \phi_{i-1}^{n+1/2} -
\alpha_i \phi_{i-2}^{n+1/2} -  h \partial_x^{\alpha_i-} \phi_i^n\right)}{2 + (1+\alpha_i) C_i^+} ,\\[1.5ex]
i=0,1,\ldots,I \,,
\end{array} \\[2ex]
\label{1d2o_kappaelegant2}
\begin{array}{l}
\phi_i^{n+1} = \dfrac{2 \phi_i^{n+1/2} - C_i^- \left((1+2 \alpha_i) \phi_{i+1}^{n+1}  - \alpha_i \phi_{i+2}^{n+1} + h \partial_x^{\alpha_i+} \phi_i^{n+1/2}\right)}{2 - (1+\alpha_i) C_i^-} ,\\[1.5ex]
i=I,I-1,\ldots,0 \,,
\end{array}
\end{eqnarray}
if the following replacements are used - instead of (\ref{1d2o_kappaelegant1}) one uses (\ref{1dbcn}) for $i=0$ and $C_0\ge 0$ and (\ref{zero}) for $i=1+i^*$ in (\ref{sign}), and, analogously, instead of (\ref{1d2o_kappaelegant2}) one uses (\ref{1dbcnB}) for $i=I$ and $C_I< 0$ and (\ref{zero}) for $i=i^*$ in (\ref{sign}). 
\end{proposition}
\begin{proof}

We rewrite (\ref{1d2obetweenelegant}) to the form 
\begin{eqnarray}
\label{1d2o_kappa}
\begin{array}{l}
\phi_i^{n+1} = \left(1 + \frac{1+\alpha}{2} |C_i| \right)^{-1}  
\Big( \phi_{i-1}^{n+1} \frac{1+2 \alpha}{2}  C_i^+ -
\phi_{i-2}^{n+1} \frac{\alpha}{2} C_i^+  \\[2ex]
- \, \phi_{i+1}^{n+1}\frac{1+2 \alpha}{2} C_i^- + \phi_{i+2}^{n+1} \frac{\alpha}{2} C_i^-  + 
\phi_i^n - C_i^+ \frac{h}{2} \partial_x^{\alpha-} \phi_i^n - C_i^- \frac{h}{2} \partial_x^{\alpha +} \phi_i^n \Big ) .
\end{array}
\end{eqnarray}

We divide the ordered set $\{0,1,\ldots,I \}$ into distinct ``uniterrupted'' ordered subsets $\mathcal{I}^+_k$ and $\mathcal{I}^-_k$ such that $C_i\ge 0$ for $i \in \mathcal{I}^+_k$ and $C_i<0$ for $i \in \mathcal{I}^-_k$ for $k=0,1,\ldots$. Let $\mathcal{I}_k^+$ be non-empty for some $k$ and $i_k$ be its first index. Then one must obtain either $i_k=0$ or $i_k=1+i^*$ in (\ref{sign}) and consequently the value $\phi_i^{n+1}$ is defined explicitly by either (\ref{1dbcn}) or (\ref{zero}). Analogously if $\mathcal{I}_k^-$ is non-empty for some $k$, then for its last index $i_k$ one must obtain either $i_k=I$ or $i_k=i^*$ in (\ref{sign}) and the value $\phi_i^{n+1}$ is defined explicitly by either (\ref{1dbcnB}) or (\ref{zero}). 

Clearly, (\ref{1d2o_kappa}) turns to (\ref{1d2o_kappaelegant1}) for $i\in\mathcal{I}^+_k$ if $i\neq i_k$ and (\ref{1d2o_kappa}) turns to (\ref{1d2o_kappaelegant2}) for 
$i\in\mathcal{I}^-_k$ if $i\neq i_k$.
Complementary, (\ref{1d2o_kappaelegant1}) for $i\in\mathcal{I}^-_k$ and (\ref{1d2o_kappaelegant2}) for $i\in\mathcal{I}^+_k$ takes the simple form $\phi_i^{n+1/2}=\phi_i^n$ and $\phi_i^{n+1}=\phi_i^{n+1/2}$, respectively. Therefore, after the forward substitution for $i=0,1,\ldots,I$ the values $\phi_i^{n+1/2}$ solve the algebraic equations (\ref{1d2o_kappa}) for $i\in\mathcal{I}_k^+$ for all existing $k$, and after the backward substitution for $i=I,I-1,\ldots,0$, the values $\phi_i^{n+1}$ solve (\ref{1d2o_kappa}) for $i\in\mathcal{I}^-_k$ and $\phi_i^{n+1}=\phi_i^{n+1/2}$ for $i\in\mathcal{I}_k^-$, so the numerical solution is completely determined.

We note that one can use $\alpha=0$ for $i=1$ if $C_1>0$ or $i=I-1$ if $C_{I-1}<0$ or some extrapolation procedures to express the values $\phi_{-1}^{n+1}$ and $\phi_{I+1}^{n+1}$, if necessary.
\end{proof}

Finally, let us briefly comment the advantages of the presented scheme with respect to analogous existing fully explicit and fully implicit schemes. Comparing to explicit ones, the presented semi-implicit scheme has no restriction on time steps due to stability as proved in Proposition \ref{propacc}. This can be used for many problems where such restriction is unpractical as discussed in Introduction. As shown in Proposition \ref{propexpl}, the forward substitution is, in fact, necessary only for non-negative Courant numbers and the backward substitution only for negative ones. Therefore, each value of numerical solution in 1D case is obtained formally using only one explicit expression. Consequently, the computational cost of the semi-implicit scheme is comparable in this case to the cost of explicit schemes with latter ones having stability restriction on discretization.

In the class of analogous parametric second order accurate (semi-) implicit schemes, the one presented here has a fully upwinded form in the implicit part for any value of the parameter. To our knowledge, this is not the case for other related schemes \cite{frolkovivc2018semi} that give systems of linear algebraic equations with matrices having less convenient properties and more involved solution procedure. This can be very convenient for nonlinear conservation laws as confirmed by first preliminary results in \cite{zeravy2021}, because the (semi-) implicit schemes can lead to non-trivially coupled nonlinear algebraic equations. Moreover, the schemes with a fixed stencil (i.e. a fixed value of $\alpha$) can lead to oscillatory numerical solution for non-smooth solutions, so variable choice of $\alpha$ can be used to suppress such unphysical oscillations \cite{zeravy2021}.

\section{Conservative advection equation}
\label{sec-c}

The linear advection equation in the conservative form is written as
\begin{equation}
    \label{1dcl}
    \partial_t \phi + \partial_x \left(v \phi \right) = 0 \,, \quad \phi(x,0) = \phi^0(x) \,.
\end{equation}
The same assumptions on the input functions, the initial and boundary conditions as in the non-conservative case apply also here.

% Other form 
% \begin{equation}
%     \label{1dcl}
%     \partial_t \phi + \partial_x f(\phi) = 0 \,, \quad \phi(x,0) = \phi^0(x) \,.
% \end{equation}
% where $f(\phi)=v \phi$.

To use a conservative finite difference (or a finite volume) method, we divide the interval $(0,L)$ to subintervals (the "control volumes") $(x_{i-1/2},x_{i+1/2})$, where the "face points" are given by  $x_{i+1/2} = i h$, $i=0,1,\ldots,I$ using the discretization step $h=L/I$. The points $x_i$ are now shifted compared to the notation in Section \ref{sec-nonc}, namely 
$$x_i = i h - h/2 \,, \,\, i=1,2,\ldots I .$$ Our aim is now to find the approximations 
$$
\Phi_i^n \approx \frac{1}{h}\int_{x_{i-1/2}}^{x_{i+1/2}} \phi(x,t^n) \, dx \approx \phi(x_i,t^n) \,,
$$
for $i=1,2,\ldots,I$. For the initial conditions, we consider $\Phi_i^0=\phi^0(x_i)$ that is a second order accurate approximation of the above integrals. The velocity is evaluated in points $x_{i-1/2}$, i.e. $v_{i+1/2} := v(x_{i+1/2})$, in particular the boundary fluxes are given by
$$
v_{1/2} = v(x_{1/2}) = v(0) \,, \,\,\, v_{I+1/2} = v(x_{I+1/2}) = v(L) \,.
$$

\subsection{First order accurate fully implicit scheme}

We define the scheme in a locally conservative form 
\begin{eqnarray}
\label{1d1ofv}
 \frac{\Phi_i^{n+1}-\Phi_i^n}{\tau} + \\[1ex]\nonumber
 \frac{v_{i+1/2}^+ \Phi_{i}^{n+1} + v_{i+1/2}^- \Phi_{i+1}^{n+1}
 - v_{i-1/2}^+ \Phi_{i-1}^{n+1} - v_{i-1/2}^- \Phi_{i}^{n+1}}{h} = 0 \,.
\end{eqnarray}
% where the values $\Phi_{i\pm 1/2}^{n+1}$ for inner faces are chosen in the upwind way, 
% \begin{eqnarray}
% \label{fluxes}
% \Phi_{i+1/2}^{n+1} = \left \{
% \begin{array}{lr}
%  \Phi_i^{n+1} \,, \,\,   &  v_{i+1/2} \ge 0 \\[1ex]
%   \Phi_{i+1}^{n+1} \,, \,\,  & v_{i+1/2} < 0
% \end{array}
% \right . \,.
% \end{eqnarray}
% or 
% \begin{eqnarray}
% \label{fluxes}
% F_{i-1/2}^{n+1} = [v_{i-1/2}]^+ \Phi_{i-1}^{n+1} +
% [v_{i-1/2}]^- \Phi_{i}^{n+1}
% \end{eqnarray}

Indexing now the signed Courant numbers at the faces
\begin{equation}
\nonumber
C_{i+1/2} := \frac{\tau v_{i+1/2}}{h} \,, \,\, i=0,1,\ldots,I \,, 
\end{equation}
we can rewrite the equations (\ref{1d1ofv}) for $i=1,2,\ldots,I$ to the form
\begin{eqnarray}
\label{1d1ofvr}
\Phi_i^{n+1} = \frac{\Phi_i^n -  C_{i+1/2}^-\Phi_{i+1}^{n+1} + C_{i-1/2}^+\Phi_{i-1}^{n+1}}{1  + C_{i+1/2}^+ - C_{i-1/2}^-}  \,,
\end{eqnarray}
where we define, formally,
\begin{eqnarray}
\nonumber
\Phi_{0}^{n+1} = 
\phi_0(t^{n+1}) , \quad
\Phi_{I+1}^{n+1} = 
\phi_L(t^{n+1}) \,.
\end{eqnarray}

The scheme (\ref{1d1ofvr}) defines the first order accurate conservative implicit upwind method for (\ref{1dcl}). The scheme represents a system of linear algebraic equations that can be solved using one forward and one backward substitution as described in the previous section. This can be viewed as the most important advantage of the first order accurate upwind scheme together with its locally conservative form. It is important to note that the case of zero velocity with diverging characteristics as described in Remark \ref{rem1} is captured by the conservative scheme (\ref{1d1ofvr}) automatically and no special treatment is required here.

The scheme can be written using two fractional time steps with the prescribed order,
\begin{eqnarray}
\nonumber
\Phi_i^{n+1/2} = \left( 1  + C_{i+1/2}^+ \right)^{-1} \left( \Phi_i^n  + C_{i-1/2}^+\Phi_{i-1}^{n+1/2} \right) \,,\,\,
i=0,1,\ldots,I
\\[2ex]
\nonumber
\Phi_i^{n+1} = \left( 1   - C_{i-1/2}^- \right)^{-1} \left(\Phi_i^{n+1/2} -  C_{i+1/2}^-\Phi_{i+1}^{n+1} \right) \,,\,\,
i=I,I-1,\ldots,0 \,.
\end{eqnarray}

The main disadvantage is again the low accuracy that motivates us to extend the scheme in a form of second order accurate semi-implicit method.

\subsection{Parametric class of second order accurate semi-implicit schemes}

To derive the error term of the first order accurate scheme (\ref{1d1ofv}), we express again the values $\phi(x_i,t^n)$ using Taylor series at $(x_i,t^{n+1})$. Instead of (\ref{taylor}), we obtain now
\begin{eqnarray}
    \nonumber
    \phi(x_i,t^n) = 
    \phi(x_i,t^{n+1}) + \\[1ex] \nonumber
    \tau \partial_x (v(x_i) \phi(x_i,t^{n+1})) - \frac{\tau^2}{2} \partial_{tx} \left(v(x_i)  \phi(x_i,t^{n+1}) \right) + \mathcal{O}(\tau^3) \,,
\end{eqnarray}
where we exploited that $\partial_t \phi=- \partial_x (v \phi)$ and $\partial_{tt} \phi=- \partial_{tx} (v \phi)$. Using additionally (\ref{taylorx}) and
$$
v_{i\pm 1/2} = v_i \pm \frac{h}{2} v'_i \,, \quad \partial_x(v_i \phi(x_i,t^{n+1})) = 
v_i \partial_x \phi(x_i,t^{n+1}) + v'_i \phi(x_i,t^{n+1})\,,
$$
we obtain the following form of the second order error term
\begin{eqnarray}
\nonumber
e_2(x_i,t^{n+1}) := - \frac{\tau^2}{2} \partial_{tx} (v_i \phi(x_i,t^{n+1})) \pm
\frac{\tau h}{2} \partial_x \left(v_i \partial_x \phi(x_i,t^{n+1}) \right)\,,
\end{eqnarray}
where $\pm = \sgn(v_i)$. Now denoting analogously to (\ref{kappa})
$$
\Phi_{i,\alpha-}^* = \alpha \Phi_i^* + (1-\alpha) \Phi_{i+1}^* \,, \quad
\Phi_{i,\alpha+}^* = (1-\alpha) \Phi_i^* + \alpha  \Phi_{i+1}^* \,, \quad
$$
we can apply the following approximations,
\begin{eqnarray}
\nonumber
\partial_{x} (v(x_i) \phi(x_i,t^{n+1})) \approx \\[1ex] \nonumber
 v_{i+1/2}^+ \, \Phi_{i,\alpha-}^{n+1} - v_{i-1/2}^+ \, \Phi_{i-1,\alpha-}^{n+1} + 
  v_{i+1/2}^- \, \Phi_{i,\alpha+}^{n+1} - v_{i-1/2}^- \, \Phi_{i-1,\alpha+}^{n+1} \,,
\end{eqnarray}
and
\begin{eqnarray}
\nonumber
h v_{i+1/2} \partial_x \phi(x_{i+1/2},t^{n+1}) \approx \\[1ex]\nonumber
v_{i+1/2}^+ ( \Phi_{i,\alpha-}^{n+1} - \Phi_{i-1,\alpha-}^{n+1}) + v_{i+1/2}^- ( \Phi_{i+1,\alpha+}^{n+1} - \Phi_{i,\alpha+}^{n+1} ) \,.
\end{eqnarray}

Together with the above approximations we use the backward finite difference in time and the central difference for the first space derivative in the second term of $e_2$. After some algebraic manipulations when several terms cancel, we obtain
\begin{eqnarray}
\nonumber
e_2(x_i,t^{n+1}) \approx \\[1ex]\nonumber
- \frac{\tau}{2 h} \left(v_{i+1/2}^+ \left(\Phi_{i-1,\alpha-}^{n+1} - \Phi_{i,\alpha-}^n \right) +
v_{i+1/2}^- \left(\Phi_{i+1,\alpha+}^{n+1} - \Phi_{i,\alpha+}^{n} \right) \right. \\[1ex] \nonumber \left.
- v_{i-1/2}^+ \left(\Phi_{i-2,\alpha-}^{n+1} - \Phi_{i-1,\alpha-}^n \right) 
- v_{i-1/2}^- \left(\Phi_{i,\alpha+}^{n+1} - \Phi_{i-1,\alpha+}^{n} \right) \right)\,.
\end{eqnarray}

Putting together the first order and the second order approximations, the second order accurate semi-implicit conservative scheme can be written in the form
\begin{equation}
\label{fv2nd}
\Phi_i^{n+1} + \frac{\tau}{h} \left( F_{i+1/2}^{n+1/2} -  F_{i-1/2}^{n+1/2} \right) = \Phi_i^n\,,
\end{equation}
where the numerical fluxes for $i=0,1,\ldots,I$ are defined by
\begin{eqnarray}
\label{F}
F_{i+1/2}^{n+1/2} := v_{i+1/2}^+ \left( \Phi_i^{n+1} - \frac{1}{2}\Phi^{n+1}_{i-1,\alpha-} + \frac{1}{2}\Phi^{n}_{i,\alpha-} \right) +\\[1ex]
\nonumber
v_{i+1/2}^- \left( \Phi_{i+1}^{n+1} - \frac{1}{2} \Phi^{n+1}_{i+1,\alpha+} +
\frac{1}{2} \Phi^{n}_{i,\alpha+} \right) \,.
\end{eqnarray}

Concerning the values $\Phi_{0}^{n+1}$ and $\Phi_{I+1}^{n+1}$ that occur in (\ref{F}) for $i=1$ if $v_{3/2}>0$, and for $i=I$ if $v_{I-1/2}<0$, we use the linear extrapolation,
\begin{equation}
    \label{extrap}
    \Phi_{0}^{n+1} = 2 \phi_0(t^{n+1}) - \Phi_1^{n+1} \,, \quad
    \Phi_{I+1}^{n+1} = 2 \phi_L(t^{n+1}) - \Phi_I^{n+1} \,.
\end{equation}
Analogously, if $v_{1/2}<0$ and $v_{I+1/2}>0$, we again linearly extrapolate the missing values by
\begin{equation}
    \nonumber
    \Phi_0^n = 2 \phi_0(t^n) - \Phi_1^n \,, \quad
    \Phi_{I+1}^n = 2 \phi_L(t^n) - \Phi_I^n \,.
\end{equation}

Concerning the inflow fluxes $F_{1/2}^{n+1/2}$ and $F_{I+1/2}^{n+1/2}$ at the boundary, we apply a second order accurate approximation, e.g.,
\begin{eqnarray*}
    F_{1/2}^{n+1/2} = v_{1/2} \phi_0(t^{n+1/2}) , \hbox{  if } v_{1/2} \ge 0 \,,  \quad \\[1ex]
    F_{I+1/2}^{n+1/2} = v_{I+1/2} \phi_L(t^{n+1/2}) , \hbox{  if } v_{I+1/2} \le 0 \,.
\end{eqnarray*}
% \begin{eqnarray}
%     \label{bf2o}
%     F_{1/2}^{n+1/2} = v_{1/2} \frac{1}{2} (\phi_0(t^{n})+\phi_0(t^{n+1})) , \hbox{  if } v_{1/2} \ge 0 \,,  \quad \\[1ex]
%     F_{I+1/2}^{n+1/2} = v_{I+1/2} \frac{1}{2} (\phi_L(t^{n})+\phi_L(t^{n+1})) , \hbox{  if } v_{I+1/2} \le 0 \,.
% \end{eqnarray}

We formulate now the scheme (\ref{fv2nd}) in the form suitable for the fast sweeping method. First we substitute (\ref{F}) to (\ref{fv2nd}),
\begin{eqnarray}
\nonumber
\Phi_i^{n+1} + \frac{1}{2} C_{i+1/2}^+ \left( (1+\alpha) \Phi_i^{n+1} - \alpha \Phi^{n+1}_{i-1}  + \alpha \Phi_i^n + (1-\alpha) \Phi_{i+1}^n \right) \\[1ex]
\nonumber
+ \left.
\frac{1}{2} C_{i+1/2}^- \left( (1+\alpha) \Phi_{i+1}^{n+1} - \alpha \Phi^{n+1}_{i+2} + (1-\alpha) \Phi_i^n + \alpha \Phi_{i+1}^n \right)  
 \right. - \\[1ex]
\nonumber
\frac{1}{2} C_{i-1/2}^+ \left( (1+\alpha) \Phi_{i-1}^{n+1} - \alpha \Phi^{n+1}_{i-2}  
  + \alpha \Phi_{i-1}^n + (1-\alpha) \Phi_{i}^n\right) -\\[1ex]
\nonumber
\frac{1}{2} C_{i-1/2}^- \left( (1+\alpha) \Phi_{i}^{n+1} - \alpha \Phi^{n+1}_{i+1}  
  + (1-\alpha) \Phi_{i-1}^n + \alpha \Phi_{i}^n\right)  = \Phi_i^n\,.
\end{eqnarray}

Collecting all terms, we can write
\begin{eqnarray}
\label{fsmfv}
\left( 2 + (1+\alpha) C_{i+1/2}^+ - (1+\alpha) C_{i-1/2}^- \right) \Phi_i^{n+1} =  2 \Phi_i^n  \\[1ex] \nonumber
+ \left.  C_{i+1/2}^+ \left( \alpha \Phi_{i-1}^{n+1} - \alpha \Phi_i^n - (1-\alpha) \Phi_{i+1}^n \right) \right. \\[1ex] \nonumber - 
\left. C_{i+1/2}^-\left((1+\alpha) \Phi_{i+1}^{n+1} - \alpha \Phi_{i+2}^{n+1}  + \alpha \Phi_{i+1}^n + (1-\alpha) \Phi_i^n \right) \right. \\[1ex]
\nonumber
+ \left. C_{i-1/2}^+\left((1+\alpha) \Phi_{i-1}^{n+1}-\alpha \Phi_{i-2}^{n+1}  + \alpha \Phi_{i-1}^n + (1-\alpha) \Phi_{i}^n \right) \right. \\[1ex] \nonumber
- C_{i-1/2}^- \left( \alpha \Phi_{i+1}^{n+1} - \alpha \Phi_{i}^n- (1-\alpha) \Phi_{i-1}^n \right)\,.
\end{eqnarray}
Dividing (\ref{fsmfv}) by the term before $\Phi_i^{n+1}$ we obtain the formula to be used with the fast sweeping method. In the case that $C_{i-1/2}<0$ and $C_{i+1/2}>0$ one can use the first order scheme (\ref{1d1ofvr}) instead of (\ref{fsmfv}).

For $i=1$ and $i=I$ we modify the scheme according to the boundary conditions as mentioned before, namely for $i=1$ 
\begin{eqnarray}
\label{fsmfv1}
\left( 2 + (1+2\alpha) C_{3/2}^+ - (1+\alpha) C_{1/2}^- \right) \Phi_1^{n+1} =  2 \Phi_1^n \\[1ex] \nonumber
+ \left.  C_{3/2}^+ \left( 2 \alpha \phi_{0}(t^{n+1}) - \alpha \Phi_1^n - (1-\alpha) \Phi_{2}^n \right) \right. \\[1ex] \nonumber - 
\left. C_{3/2}^-\left((1+\alpha) \Phi_{2}^{n+1} - \alpha \Phi_{3}^{n+1}  + \alpha \Phi_{2}^n + (1-\alpha) \Phi_1^n \right)   \right. \\[1ex] \nonumber
+ \left. C_{1/2}^+\phi_{0}(t^{n+1/2}) - C_{1/2}^- \left( \alpha \Phi_{2}^{n+1} + (1-2\alpha) \Phi_{1}^n -2 (1- \alpha) \phi_{0}(t^{n}) \right) \right.\,,
\end{eqnarray}
and analogously for $i=I$,
\begin{eqnarray}
\label{fsmfvI}
\left( 2 + (1+\alpha) C_{I+1/2}^+ - (1+2\alpha) C_{I-1/2}^- \right) \Phi_I^{n+1} =  2 \Phi_I^n \quad \\[1ex] \nonumber
+ \left.  C_{I+1/2}^+ \left( \alpha \Phi_{I-1}^{n+1} + (1-2\alpha) \Phi_I^n - 2 (1-\alpha) \phi_L(t^{n}) \right) -  C_{I+1/2}^- \phi_L(t^{n+1/2})  \right. \\[1ex] \nonumber 
+ \left. C_{I-1/2}^+\left((1+\alpha) \Phi_{I-1}^{n+1}-\alpha \Phi_{I-2}^{n+1}  + \alpha \Phi_{I-1}^n + (1-\alpha) \Phi_{I}^n \right) \right. \\[1ex] \nonumber
- C_{I-1/2}^- \left( 2 \alpha \phi_L(t^{n+1}) - \alpha \Phi_{I}^n- (1-\alpha) \Phi_{I-1}^n \right)\,.
\end{eqnarray}
Note that for $i=2$ with $C_{3/2}>0$ and $i=I-1$ with $C_{I-1/2}<0$ one has to use also (\ref{extrap}).

The derived scheme is exactly mass conservative at the discrete level due to (\ref{fv2nd}) and it is second order accurate as it is suggested also by numerical experiments. In the case of constant velocity $v$ in (\ref{1dcl}) the scheme is equivalent to the parametric family of nonconservative semi-implicit scheme (\ref{1d2o_kappa}) when also the von Neumann unconditional stability is valid.

We note that although the parameter $\alpha$ can be chosen freely in each time interval, e.g. $\alpha \in [0,1]$, it shall not vary with respect to $i$ if the discrete form of local mass balance property shall be fulfilled. 

\section{Two dimensional case}
\label{sec-2d}

We now apply the so called Strang splitting to solve the non-conservative advection equation in two-dimensional case, but the idea can be applied in more dimensional cases and for the conservative form, too. 

The advection equation is now given in the form
\begin{equation}
    \label{2dadvection}
    \partial_t \phi + \vec{v} \cdot \nabla \phi = 0 \,, \quad
    \phi(x,y,0) = \phi^0(x,y) \,, (x,y) \in \Omega \,,
\end{equation}
where $\Omega \subset R^2$ takes here a simple form of a square and $\vec{v} = (v_1(x,y),v_2(x,y))$ is a given velocity vector field. The boundary conditions are defined depending on the flow regime at the boundary $\partial \Omega$ with the values prescribed only at the inflow part,
\begin{equation}
    \label{2dinflow}
    \phi(x,y,t) = \phi_{in}(x,y,t) \hbox{   if   } \vec{n}(x,y) \cdot \vec{v}(x,y) \le 0 \,, \,\,
    (x,y) \in \partial \Omega \,,
\end{equation}
where $\vec{n}$ is the outward normal vector.

We use analogous notation to derive numerical approximations as in the previous section with the addition that $y_j = j h$ for $j=0,1,\ldots,I$. Furthermore, we denote $t^{n+1/2} = t^n + \tau/2$.

The idea of the time splitting method is to approximate and split the problem (\ref{2dadvection}) into two subproblems that are coupled only by the choice of initial conditions, and that are solved separately in a specified sequence.  Let us explain the simplest variant in details.

The first subproblem takes the form of one dimensional advection equations for the parameter $y \in (0,L)$,
\begin{eqnarray}
    \label{2dsplitx}
    \partial_t \phi(x,y,t) + v_1(x,y) \partial_x \phi(x,y,t) = 0 \,, \,\, x \in (0,L) \,,
\end{eqnarray}
and the second subproblem takes the analogous form for the parameter $x \in (0,L)$,
\begin{eqnarray}
    \label{2dsplity}
    \partial_t \phi(x,y,t) + v_2(x,y) \partial_y \phi(x,y,t) = 0 \,, \,\, y \in (0,L) \,.
\end{eqnarray}

The splitting in time is realized as follows. 
Let the solution $\phi(x,y,t^n)$ (or its approximation) of (\ref{2dadvection}) be available at some time $t=t^n$, $0\le n<N$. To obtain an approximation of $\phi(x,y,t^{n+1})$ we do three steps. First, the subproblem (\ref{2dsplitx}) is solved for $t\ \in (t^n,t^{n+1/2})$ with the initial condition defined by $\phi(x,y,t^n)$. Afterwards, the subproblem (\ref{2dsplity}) is solved for $t\in (t^{n},t^{n+1})$ with the initial condition defined by the solution of (\ref{2dsplitx}) at $t=t^{n+1/2}$. Finally, the subproblem (\ref{2dsplitx}) is solved now for $t \in (t^{n+1/2},t^{n+1})$ with the initial condition defined by the solution of (\ref{2dsplity}) at $t=t^{n+1}$. The result is the desired approximation of $\phi(x,y,t^{n+1})$.

To discretize the advection equation (\ref{2dadvection}) also in space, we consider the first subproblem (\ref{2dsplitx}) only for $y=y_j$, $j=0,1,\ldots,J$, and the second one (\ref{2dsplity}) only for $x=x_i$, $i=0,1,\ldots,I$. For the resulting one-dimensional advection problems we can apply the numerical method from the previous section, when each of the resulting discrete algebraic systems can be solved in one forward and one backward substitution.

\section{Numerical experiments}
\label{sec-ne}

In following numerical experiments we want to illustrate the properties of the derived semi-implicit schemes using some standard benchmarks. For the non-conservative advection, the methods are implemented in $C$-language for two-dimensional problems. The numerical solutions of one-dimensional conservative advection is implemented in Python. The experiments with optimized choice of parameters using automatic differentiation is realized with Python and its library PyTorch \cite{NEURIPS2019_9015}.

For examples having available exact solution on a whole time interval we compute the global discrete $L_1$ errors $E$ that takes in one-dimensional case the form
\begin{equation}
    \label{globerror}
    E = E(h,\tau) = h \tau \sum_{i,n} |\phi_{i}^n - \phi(x_i,t^n)| 
\end{equation}
and analogously for two-dimensional case. If an exact solution is available only at the final time $t^N=T$, we compute the error in the form
\begin{equation}
    \label{localerror}
    E^N = E^N(h,\tau) = h \sum_{i} |\phi_{i}^N - \phi(x_i,T)| 
\end{equation}
and analogously in two-dimensional case.
\subsection{One dimensional nonconservative advection}
\label{ex01}

\begin{figure}[ht]
\centering
 \includegraphics[width=0.485\linewidth]{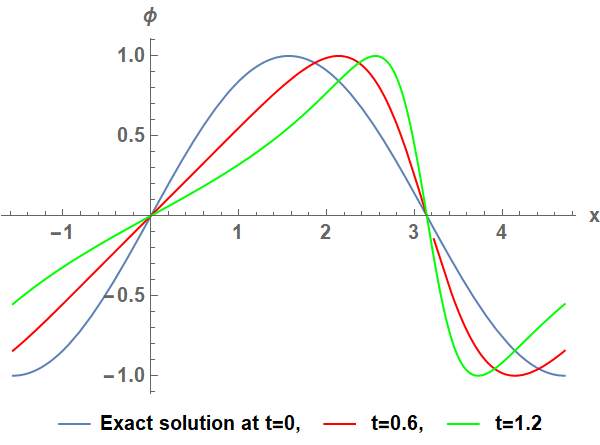}~~~~
 \includegraphics[width=0.485\linewidth]{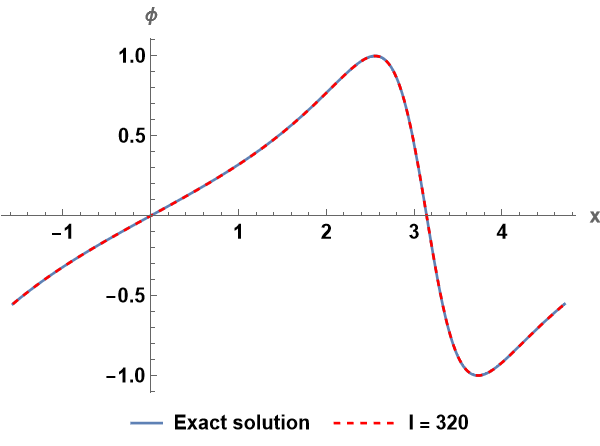}
 \caption{The exact solution for $t=0$, $t=T/2$, and $t=T$ (left) and the comparison of exact and numerical solution at $t=T$ (right) for the example \ref{ex01}.}
 \label{1d_initial}
\end{figure}

The following example contains all important features of the non-conservative linear advection equation (\ref{1dadv}) with variable velocity and general boundary conditions (\ref{1dbc}). The velocity changes its sign twice inside of computational interval with converging and expanding characteristics and it prescribes variable inflow and outflow boundary conditions. 

The example is formally treated as two-dimensional. The computational domain is $\Omega = (-\pi/2, 3\pi/2)^2$, the initial function $\phi(x,y,0)=\sin(x)$ and $T=1.2$. The velocity function $\vec{v}=(v_1(x,y),v_2(x,y))$ is depending on the spatial coordinate $x$ only with  $v_1(x,y)=\sin(x)$ and $v_2(x,y)=0$. The exact solution is given by
\begin{equation}
    \nonumber
    \phi = \sin(2 \arctan(e^{-t} \tan(x/2))) \,.
\end{equation}

We solve the example with discretization steps resulting in the maximal Courant number being approximately $3.81$. 
In Table \ref{tab1d} we present the global discrete errors $E$ in (\ref{globerror}) for two interesting choices of the parameter $\alpha$, namely $\alpha=0.5$ and $\alpha$ defined by (\ref{kappav}). One can see that the results are of a good accuracy even for Courant number larger than $1$. The EOC is approaching $2$ from above in both cases with the results for the variable $\alpha$ slightly better than for the fixed value. 

\begin{table}[ht]
	\begin{center}
	\begin{tabular}{ c c c c c c }
	\hline
 $I$ & $N$ & $E,\alpha=0.5$ & EOC & $E,\alpha$ in (\ref{kappav}) & EOC \\ 
	\hline
 40 & 1 & 0.810861 & - & 0.556925 & - \\
 80 & 2 & 0.167179 & 2.278 & 0.099711 &  2.481 \\
 160 & 4 & 0.035211 & 2.247 & 0.018519 & 2.428 \\
 320 & 8 & 0.007858 & 2.163 & 0.003831 & 2.273 \\
	\hline 
	\end{tabular}
	\end{center}
	\caption{Numerical errors and EOC for two choice of $\alpha$ in example \ref{ex01}.}
	\label{tab1d}
\end{table}

\begin{figure}[ht]
\centering
 \includegraphics[width=0.49\linewidth]{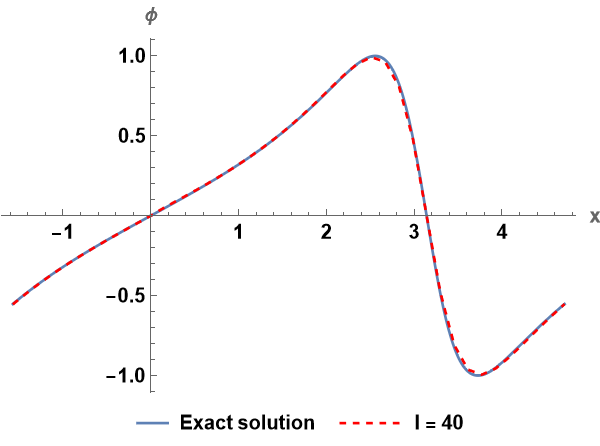}~
 \includegraphics[width=0.49\linewidth]{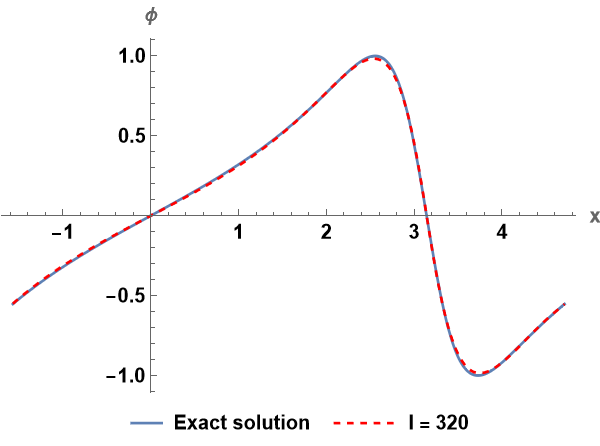}
 \caption{Comparison of exact and numerical solutions obtained only with one time step for example \ref{ex01}. The maximal Courant number is $3.81$ (left) and $30.5$ (right).}
 \label{1d_compare}
\end{figure}

Finally, we compute the example on the finest mesh using only one time step resulting in the maximal Courant number being approximately $30.5$. The obtained numerical solution is compared with the one obtained on the coarsest mesh in Figure \ref{1d_compare}. One can see that the large Courant numbers do not result in any instabilities.

\subsection{Optimization of parameters by automatic differentiation}
\label{ex01c}

To illustrate the possibilities of variable parameter $\alpha=\alpha_i^n$ in (\ref{1d2obetweenelegant}), we compute an example with variable velocity, where we let the code to optimize the values of $\alpha_i^n$ at each time step and at each grid point. To do so we use the library PyTorch \cite{NEURIPS2019_9015} to implement a straightforward gradient descent method with no constraints to minimize the loss function
$J=J(I,N) = h \tau \sum_{i,n}  (\min\{0,\phi_i^n\})^2$
with respect to parameters $\alpha_i^n$. 

In particular, we set initially $\alpha_i^n=0.5$ for $i=1,2,\ldots,I-1$ and $n=1,2,\ldots,N-1$ for which we compute the values $\phi_i^n$ of numerical solution. Afterwards using the automatic differentiation available in PyTorch we obtain the gradient of $J$ with respect to all values of $\alpha_i^n$. Next we subtract the gradient multiplied by a parameter (the "learning rate") $\eta$ from the values of $\alpha_i^n$ that results in a smaller value of $J$. In theory, one can continue with this procedure up to a point when the decrease of $J$ is not substantial. In our case, we use only one step of such optimization.

To show clearly this idea we choose the following example. The domain is $\Omega = (-2, 12)$,  $T=2\pi/\sqrt{3}$, the initial function $\phi(x,0)=\exp(-2 x^2)$ and the velocity function $v(x)=2+\sin(x)$. One can show that $\phi(x,T)=\phi^0(x-2\pi)$.

In Table \ref{tabml} we summarize the results. For three consecutively refined meshes we present the value of $J$ before and after the optimization step, and analogously the errors $E^N$ in (\ref{localerror}). One can see that the unphysical oscillations can be decreased significantly with a slight improvement in the precision.

\begin{table}[ht]
	\begin{center}
	\begin{tabular}{ r r l l l l l }
	\hline
 $I$ & $N$ & $\eta \cdot 10^6$ & $\,\, J_b \cdot 10^{-3}$ & $\,\, J_a \cdot 10^{-3}$ & $\,\, E_b^N$ & $\,\, E_a^N$ \\ 
	\hline
 70 & 50 & 0.2 & 3.68 & 0.0768 & 0.521 & 0.511 \\
 140 & 100 & 4.0 & 1.12 & 0.0156 & 0.197 & 0.190 \\
 280 & 200 & 160.0 & 0.0664 & 0.00354 & 0.0533 & 0.0448 \\
	\hline 
	\end{tabular}
	\end{center}
	\caption{The learning rate $\eta$, the values of the loss function $J$, and the errors $E^N$ in (\ref{localerror}) using three meshes for example \ref{ex01c}. The indices $b$ and $a$ means "before" and "after" optimization.}
	\label{tabml}
\end{table}

\begin{figure}[ht]
\centering
 \includegraphics[width=0.51\linewidth]{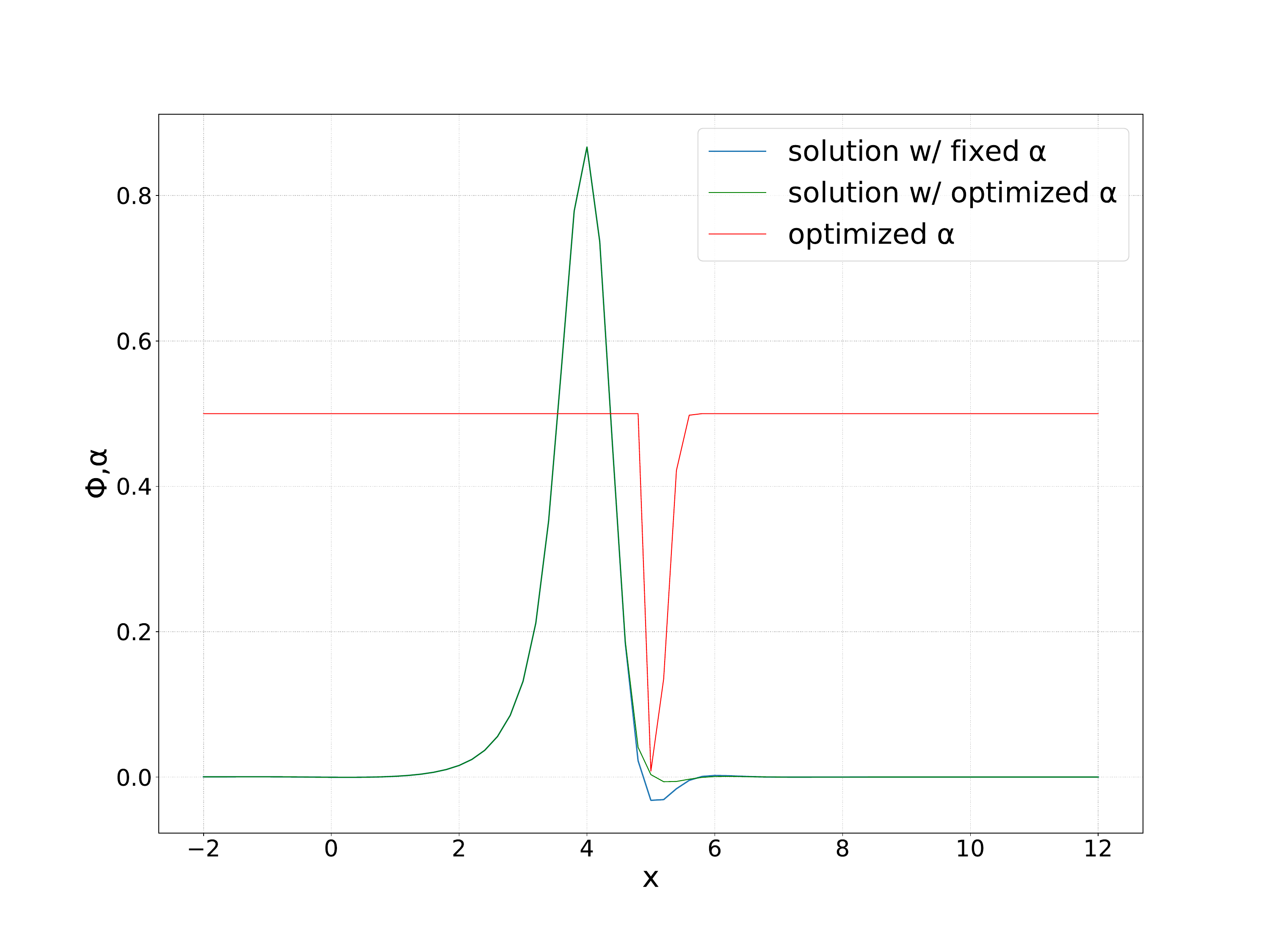}~
 \hspace{-0.06\linewidth}\includegraphics[width=0.51\linewidth]{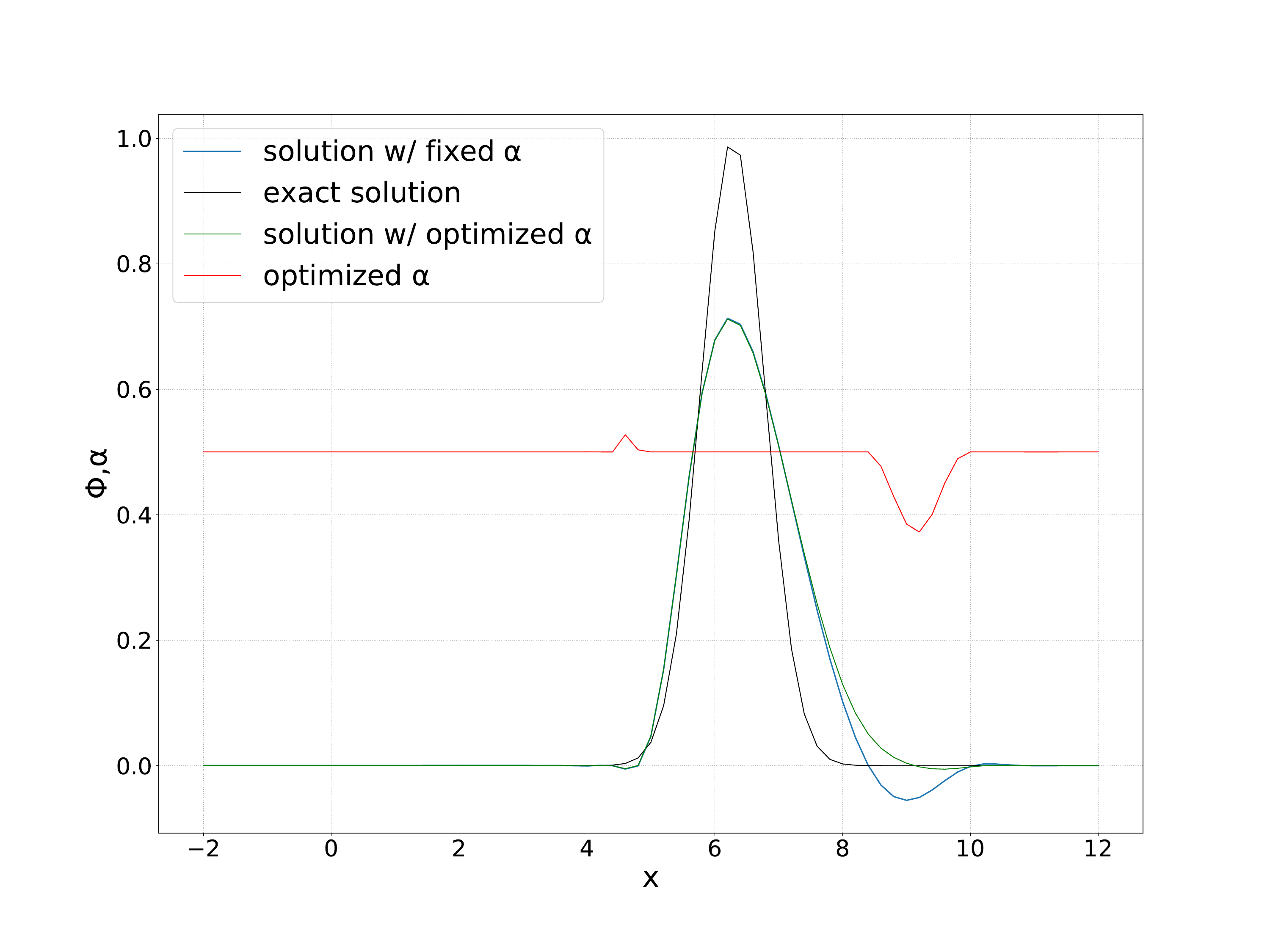}\\[-4ex]
  \includegraphics[width=0.51\linewidth]{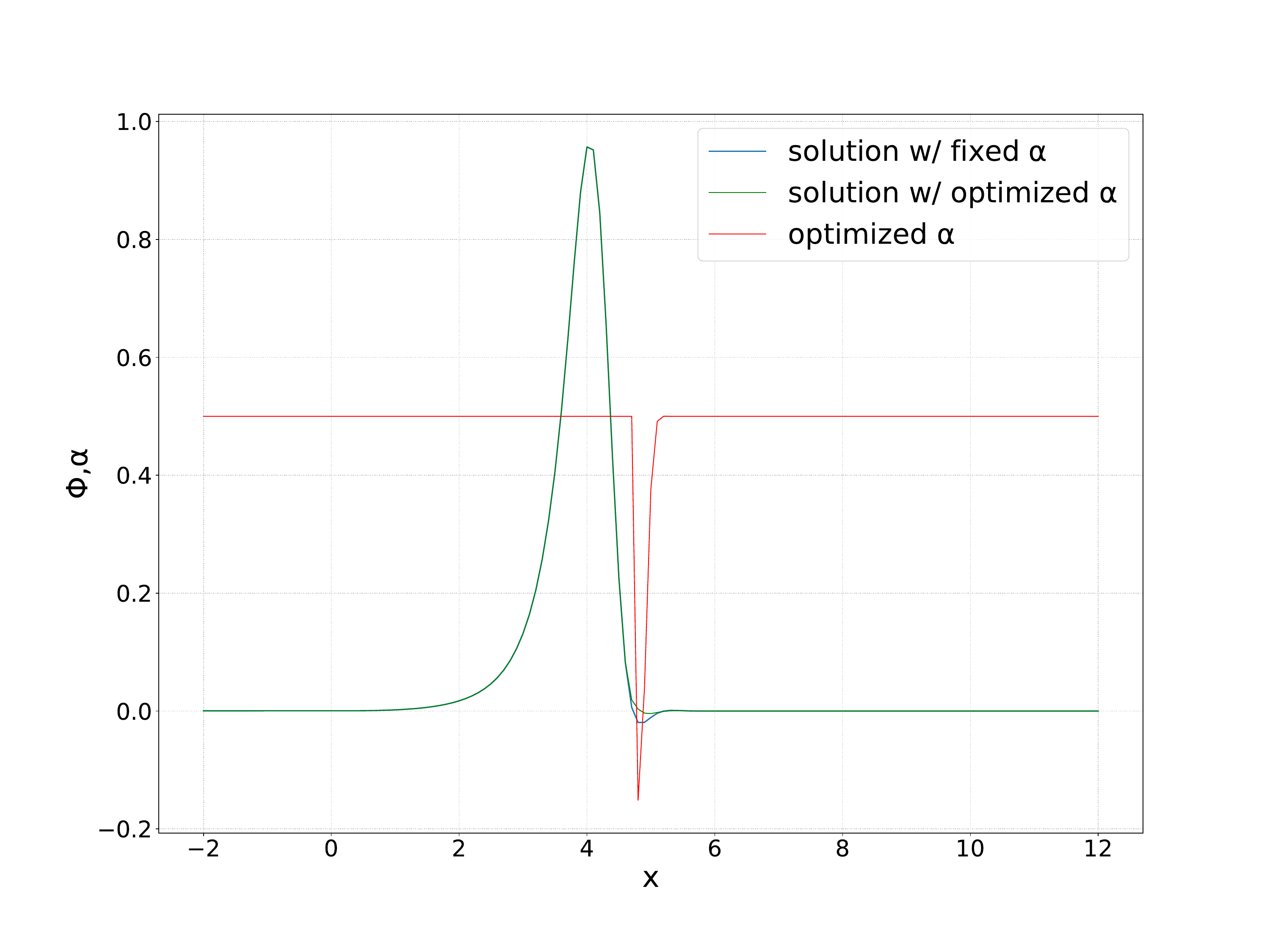}~
 \hspace{-0.06\linewidth}\includegraphics[width=0.51\linewidth]{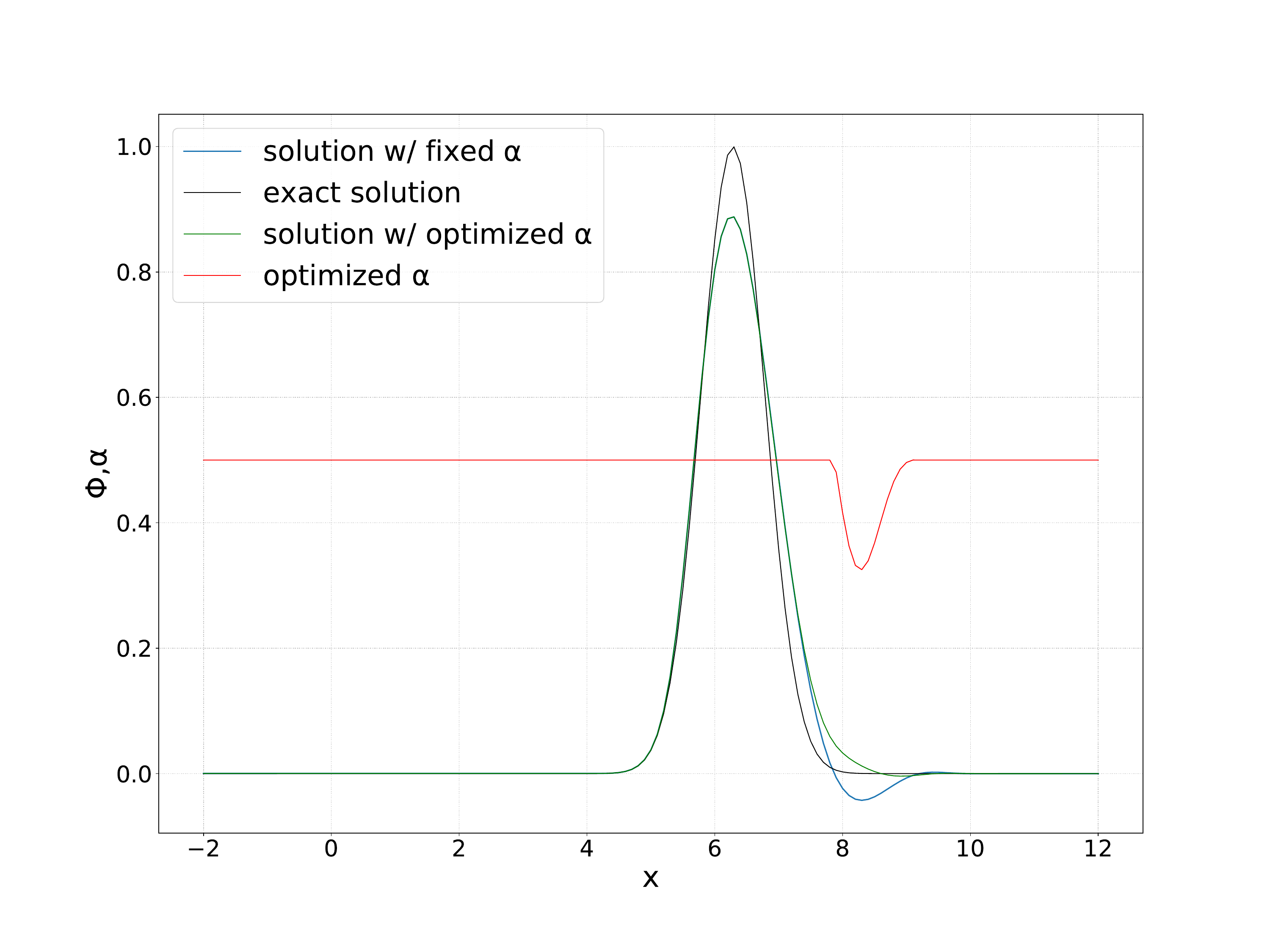}
 \caption{Numerical solution and variable parameters for $t=T/2$ (left) and for $t=T$ (right) for example \ref{ex01c}. The top row is obtained with $I=70$ and the bottom one with $I=140$.\\[-10ex]}
 \label{1d_comparec}
\end{figure}

To show the influence of optimized values $\alpha$ visually, we present the figures for $I=70$ and $N=50$ representing a rather coarse time and space discretization. In Figure \ref{1d_comparec} we present the numerical solutions at $t^{N/2}=T/2$ and $t^N=T$ for the fixed choice of all $\alpha_i^n \equiv 0.5$ when one can observe clearly some unphysical negative values. Furthermore, the numerical solutions at the same times obtained after one optimization step are plotted together with the values of $\alpha_i^{N/2}$ and $\alpha_i^N$. One can clearly observe that the largest unphysical oscillation are suppressed. Analogously, the same results are presented also for the refined discretization steps with $I=140$ and $N=100$. Note that in the latter case one can observe that $\alpha_i^n$ attained a small negative value in one grid point as we used an unconstrained minimization method. Nevertheless it caused no instabilities in the results because of appropriate choice of the loss function $J(I,N)$. The stability of our numerical scheme is proved only for non-negative values of $\alpha$, therefore, in general, one shall use a constraint optimization method especially if different type of loss function is used than in this example.\\[-6ex]

\subsection{One dimensional conservative advection}
\label{ex01b}

\begin{figure}[ht]
\centering
 \includegraphics[width=0.55\linewidth]{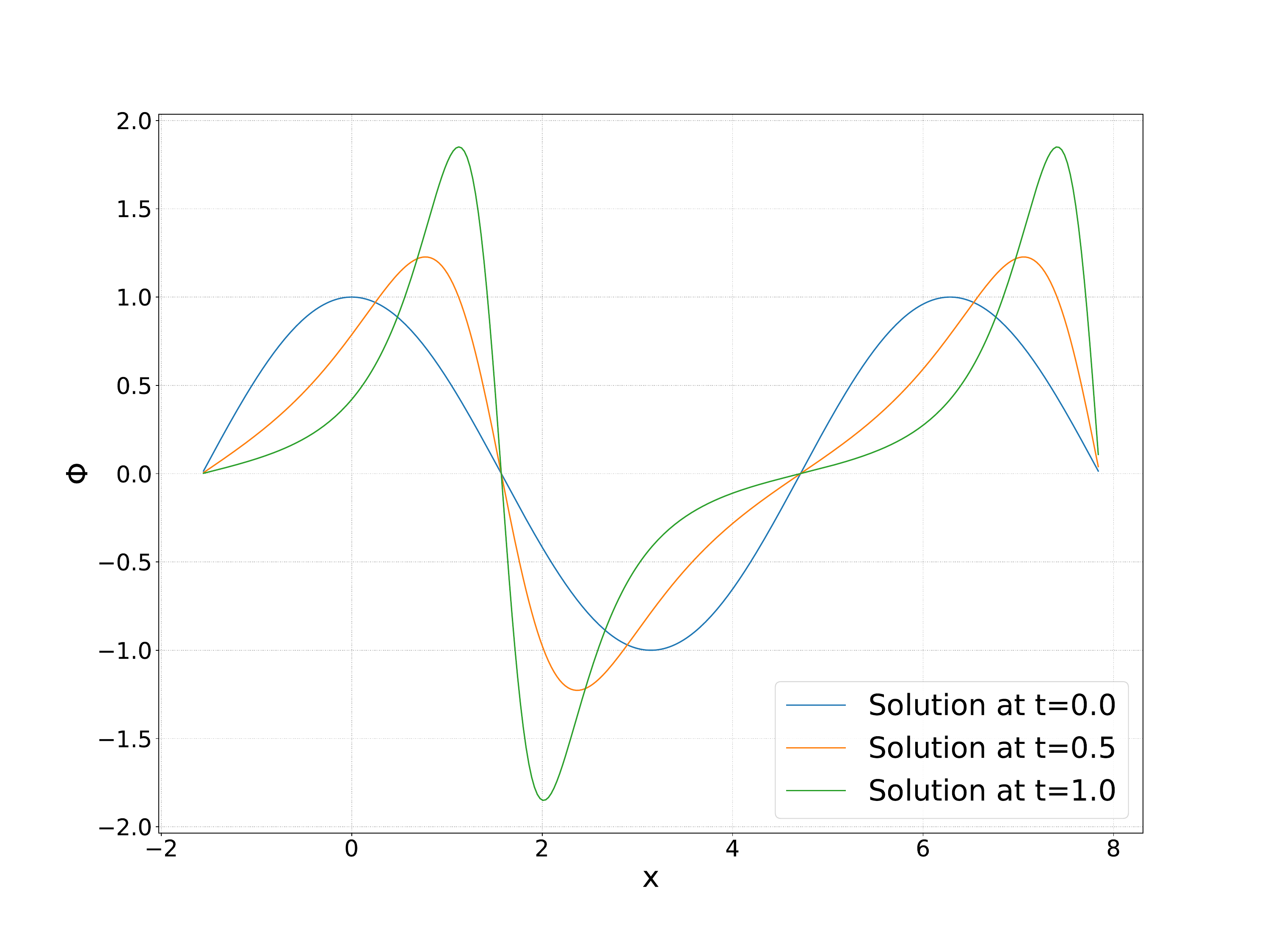}~~~~
 \hspace{-0.06\linewidth}\includegraphics[width=0.55\linewidth]{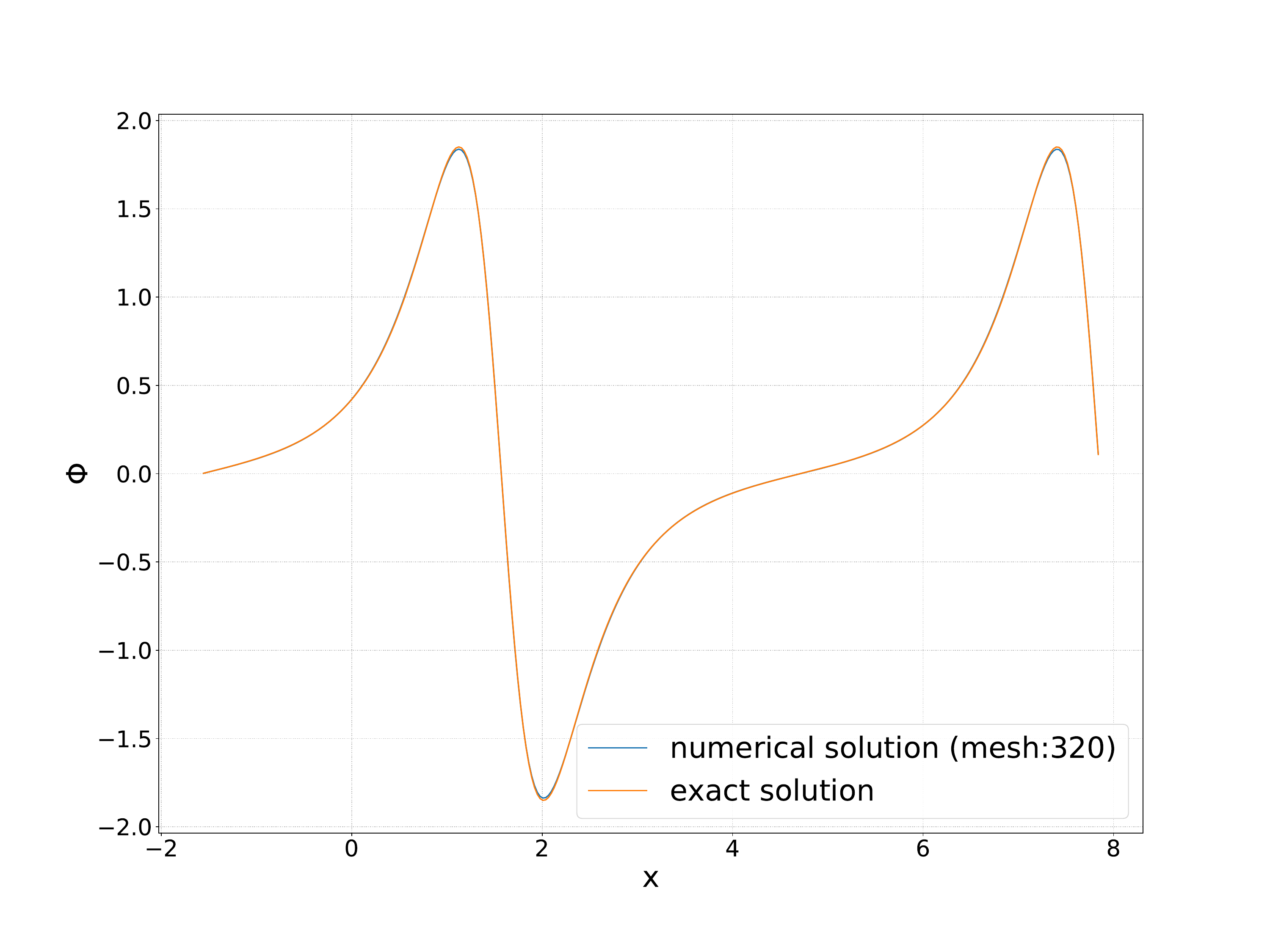}~~~~
 \caption{The exact solution at $t=0$, $t=T/2$, and $t=T$ (left) and the comparison of exact and numerical solution at $t=T$ (right) for example \ref{ex01b}.}
 \label{1d_initialb}
\end{figure}

The following example illustrates the applicability of semi-implicit scheme for the advection equation in the conservative form (\ref{1dcl}). To check the mass conservation property, we choose an example with zero velocity at boundary points. Inside of the computational interval $\Omega = (-\frac{\pi}{2}, \frac{5\pi}{2})$ the variable velocity $v(x)=\cos(x)$ changes twice its sign. The initial condition has the form $\phi(x,0)=\cos(x)$ and $T=1$. The exact solution is given by
\begin{equation}
    \nonumber
    \phi = \cos ^2\left(2 \tan ^{-1}\left(\tanh \left(\frac{1}{2} \left(t-2 \tanh ^{-1}\left(\tan
   \left(\frac{x}{2}\right)\right)\right)\right)\right)\right) \,,
\end{equation}
see the left picture in Figure \ref{1d_initialb}.\\[-6ex]

We solve the example with discretization steps such that the maximal Courant number is approximately $4.24$. The global discrete errors $E$ in (\ref{globerror}) are presented in  Table \ref{tab1db} for $\alpha=0.5$ and $\alpha=1$. The numerical results are stable even for Courant number larger than $1$, see the right picture in Figure \ref{1d_initialb} for a visual comparison. The EOC is approximately $2$ in both cases with the choice $\alpha=1$, as expected, being slightly better for the example with larger Courant numbers.

\begin{table}[ht]
	\begin{center}
	\begin{tabular}{ c c c c c c }
	\hline
 $I$ & $N$ & $E,{\alpha=0.5}$ & EOC & $E,{\alpha=1}$ & EOC \\ 
	\hline
 40 & 1 & 0.9610 & - & 0.7013 & - \\
 80 & 2 & 0.2750 & 1.81 & 0.1941 & 1.85  \\
 160 & 4 & 0.0651 & 2.08 & 0.0442 &  2.13 \\
 320 & 8 & 0.0150 & 2.12 & 0.0098 & 2.17 \\
	\hline 
	\end{tabular}
	\end{center}
	\caption{Numerical errors for example \ref{ex01b} with the maximal Courant number $4$.}
	\label{tab1db}
\end{table}

In Table \ref{tab1dc} we present the results for the same example with the four times smaller time step resulting in the maximal Courant number being approximately $1$. Of course, the precision of results increases, and, moreover, the choice $\alpha=0.5$ gives now better results than the choice $\alpha=1$.

\begin{table}[ht]
	\begin{center}
	\begin{tabular}{ c c c c c c }
	\hline
 $I$ & $N$ & $E,{\alpha=0.5}$ & EOC & $E,{\alpha=1}$ & EOC \\ 
	\hline
 40 & 4 & 0.1181 & - & 0.1683 & - \\
 80 & 8 & 0.0256 & 2.20 & 0.0461 & 1.87  \\
 160 & 16 & 0.0054 & 2.26 & 0.011 &  2.02 \\
 320 & 32 & 0.0012 & 2.16 & 0.0028 & 2.02 \\
	\hline 
	\end{tabular}
	\end{center}
	\caption{Numerical errors for example \ref{ex01b} with the maximal Courant number $1$.}
	\label{tab1dc}
\end{table}

Note that we always obtain a perfect mass conservation at the discrete level by checking
$$
\sum_i \Phi_i^n = \sum \Phi_i^0 \,, \quad n=1,2,\ldots,N 
$$
with a difference given only by rounding errors around $10^{-15}$.

\begin{figure}[ht]
\centering
 \includegraphics[width=0.55\linewidth]{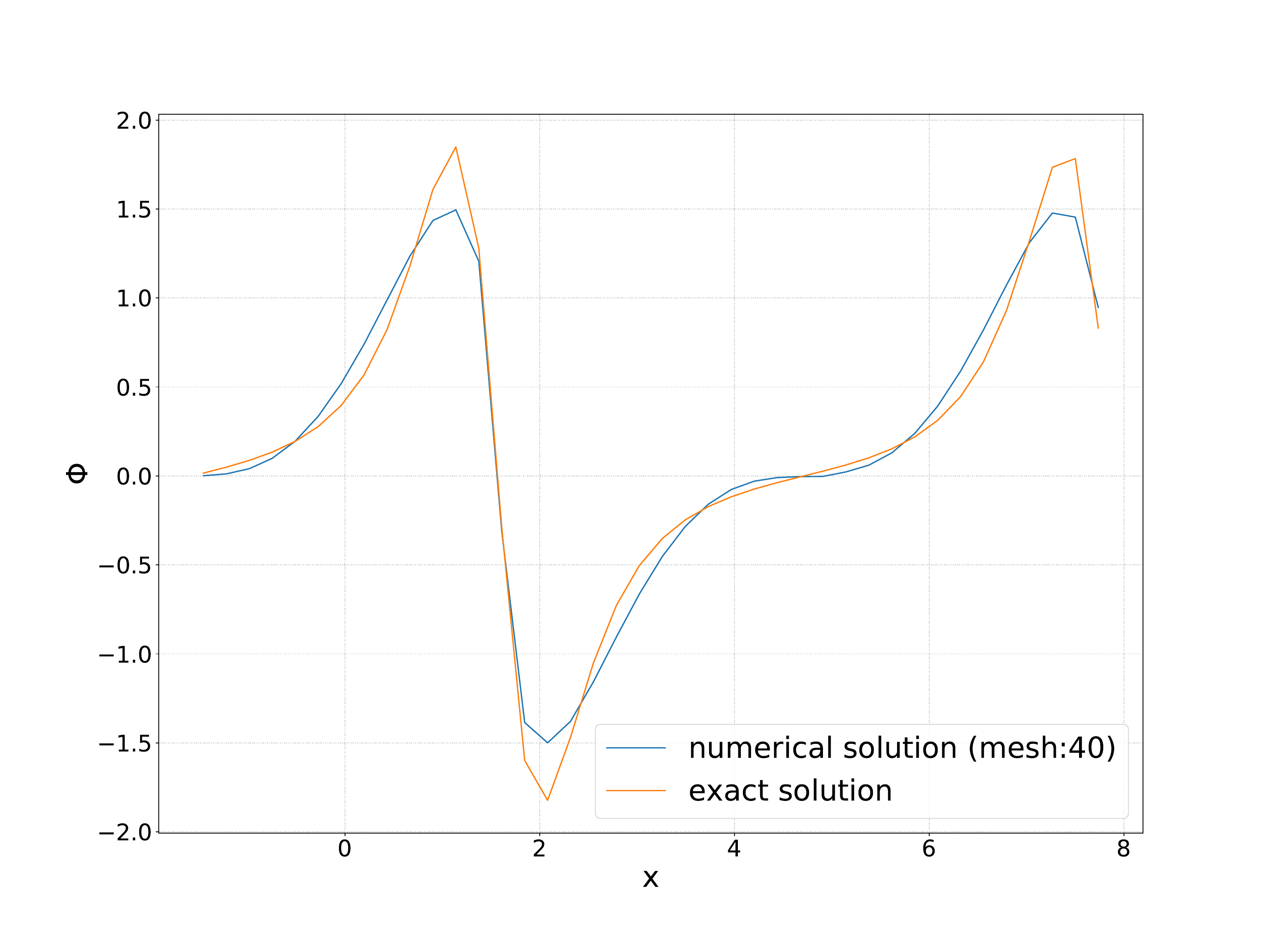}~
 \hspace{-0.06\linewidth}\includegraphics[width=0.55\linewidth]{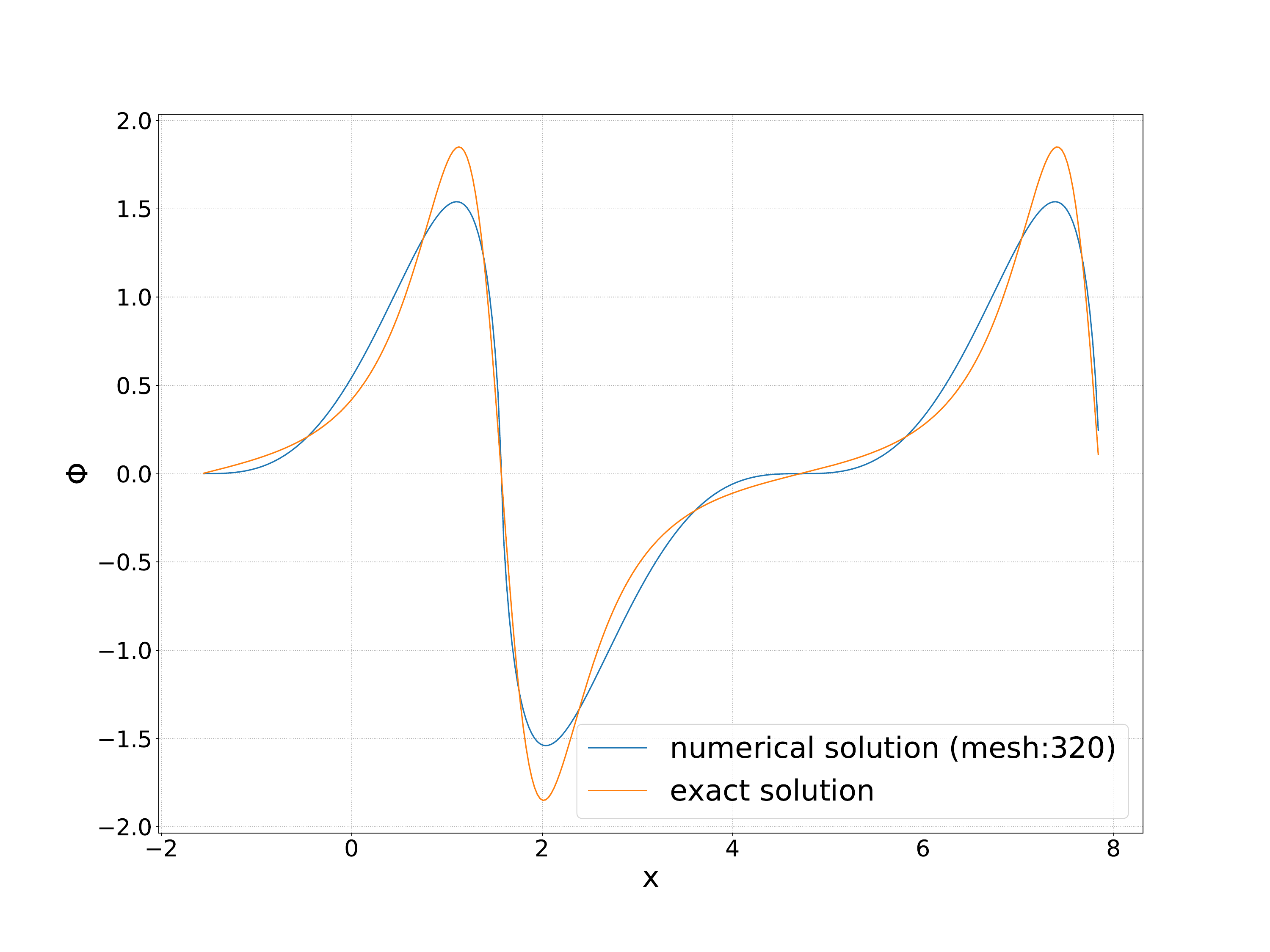}
 \caption{Numerical solutions obtained with only one time step for example \ref{ex01b}. The maximal Courant number equals $4.24$ (left) and $33.95$ (right).\\[-8ex]}
 \label{1d_compareb}
\end{figure}

Finally, to illustrate the stability of our scheme, we compute the example on the finest mesh with only one time step resulting in the maximal Courant number approximately $34$, see the right picture in Figure \ref{1d_compareb} with no instabilities occurring.\\[-6ex]

\subsection{Two-dimensional examples}

With the next example we solve an analogous problem to Example \ref{ex01} with the velocity $\vec{v}=(\sin(2\pi (x+y)/2),\sin(2\pi (x+y)/2))$ that is variable only in a diagonal direction of the  domain $\Omega = (-1,2) \times (-1,2)$. We use the initial condition $\phi(x,y,0)=\sin(2\pi (x+y)/2)$ and $T=0.24$. 

The first example is chosen to test the accuracy of the operator splitting method as described in Section \ref{sec-2d}. To compute the global error $E$ analogous to (\ref{globerror}) we use always a half time step twice in the subproblem (\ref{2dsplitx}) when compared to the time step of subproblem (\ref{2dsplity}) as described in Section \ref{sec-2d}. Consequently, the maximal Courant number in the $x$ direction is a half of the maximal Courant number in $y$ direction. Note that for explicit methods a typical stability restriction is given by the sum of directional Courant numbers that is required to be smaller than one \cite{leveque2002finite}.
We compute the example with the maximal Courant number in $y$ direction being $1.6$.  

\begin{figure}[ht]
\centering
 \includegraphics[width=0.49\linewidth]{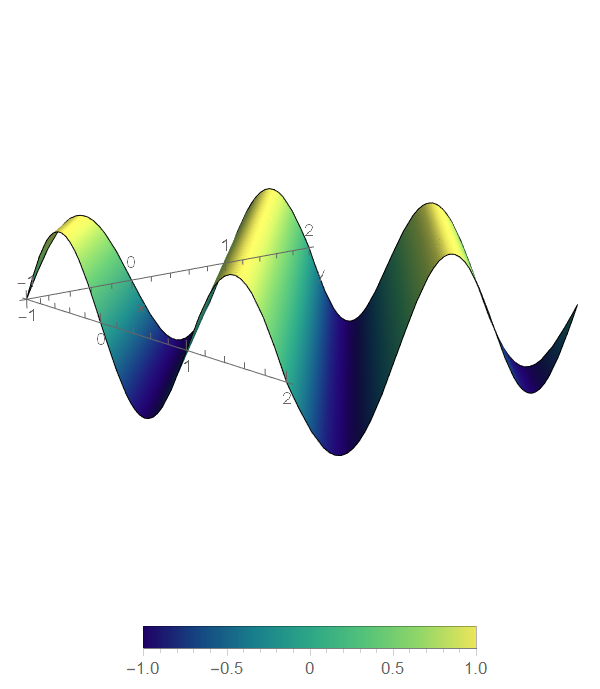}~~
 \includegraphics[width=0.49\linewidth]{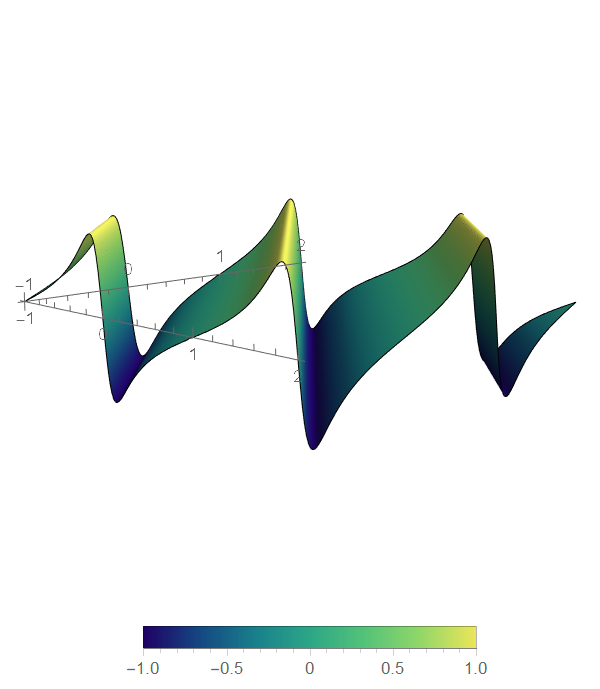}
 \caption{The initial condition for Example \ref{ex2d1} viewed from a side (left) and the corresponding numerical solution in time $t=T$ for $I=160$.}
 \label{diag_initial}
\end{figure}

The exact solution is given by
\begin{equation}
\label{ex2d1}
    \phi =\sin(2 \arctan(e^{-2\pi t} \tan(\pi \frac{x+y}{2}))).
\end{equation}
The numerical results are presented in Figure \ref{diag_initial} and in Table \ref{tabdiag}. We see that the method for two typical choices of parameter $\alpha$ is second order accurate for this example.\\[-6ex]

% \begin{figure}[ht]
% \centering
%  \includegraphics[width=0.49\linewidth]{img/test12_final_161_top.png}~~
%  \includegraphics[width=0.49\linewidth]{img/test12_final_161_side.png}
%  \caption{The numerical solution for Example \ref{ex2d1} in time $t=T$ for $I=160$ viewed from a top and a side.}
%  \label{diag_final}
% \end{figure}

\begin{table}[ht]
	\begin{center}
	\begin{tabular}{ c c l l l l }
	\hline
 $I$ & $N$ & $E,\alpha=0$ & EOC & $E,\alpha$ in (\ref{kappav}) & EOC \\ 
	\hline
 20 & 1 & 0.0874 & - & 0.0838 & - \\
 40 & 2 & 0.0179 & 2.29 & 0.0173 & 2.27 \\
 80 & 4 & 0.00319 & 2.49 & 0.00302 & 2.52 \\
 160 & 8 & 0.000624 & 2.36 & 0.000569 & 2.41 \\
	\hline 
	\end{tabular}
	\end{center}
	\caption{Numerical errors for two choices of $\alpha$ for example (\ref{ex2d1}).}
	\label{tabdiag}
\end{table}

% \subsection{Example with deformation velocity }
% \label{ex04}

Next, to test our method for a nontrivial case, we choose two-dimensional example with a deformation velocity in which the initial profile of solution is deformed significantly in time. For this example we check quantitatively not only the numerical errors, but also two other numerical artifacts - negative unphysical oscillations and a violation of mass conservation in a discrete form. As we show, the both of them are visible for a coarse mesh, but these numerical errors decrease rapidly with the mesh refinement. 

\begin{figure}[ht]
\centering
 \includegraphics[width=0.6\linewidth]{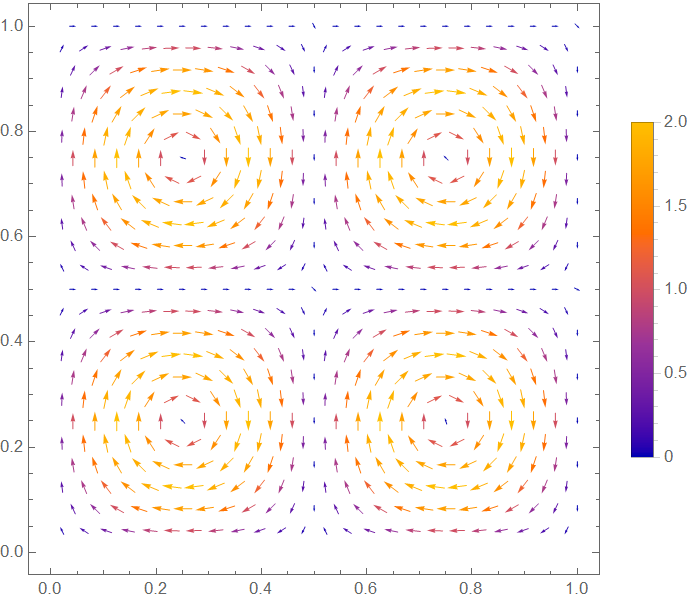}
 \caption{The velocity field $\vec{v}$ at $t=0$ for example (\ref{rotvelocity}).}
 \label{quadvelocity}
\end{figure}

The domain is the unit square $\Omega = (0,1) \times (0,1)$ and $T=1$. The velocity vector $\vec{v}$, see Figure \ref{quadvelocity}, is defined by
\begin{equation}
\label{rotvelocity}
\begin{array}{c}
  v_1(x,y,t) = -4 \cos(\pi t) \sin^2(2 \pi x) \sin(2\pi y) \cos(2 \pi y) \\[1ex]
  v_2(x,y,t) = 4 \cos(\pi t) \sin^2(2 \pi y) \sin(2 \pi x) \cos(2 \pi x)
\end{array}
\end{equation}
Note that the time dependency of $\vec{v}$ in numerical simulations is resolved by evaluating $v$ for each $n$ at $t=t^n+\tau/2$. The velocity has the zero divergence and it is equal zero at the boundary, so the integral of the initial function (the mass) shall be conserved in time.

\begin{figure}[ht]
\centering
 \includegraphics[width=0.49\linewidth]{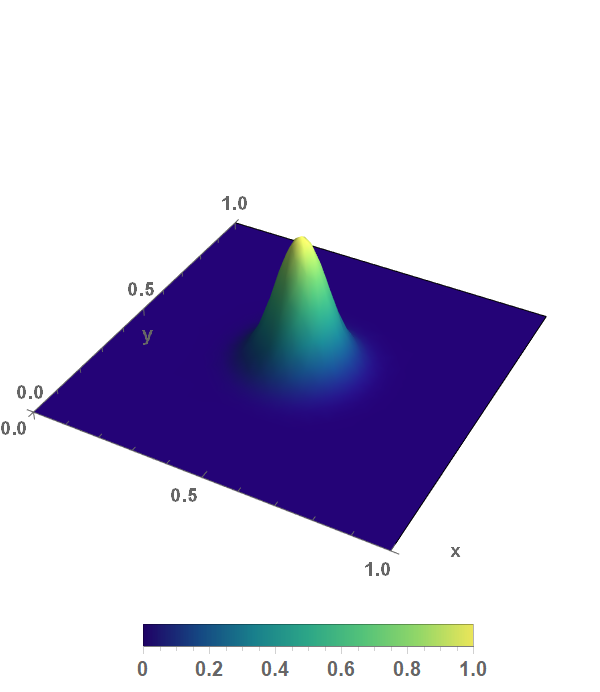}~~
 \includegraphics[width=0.49\linewidth]{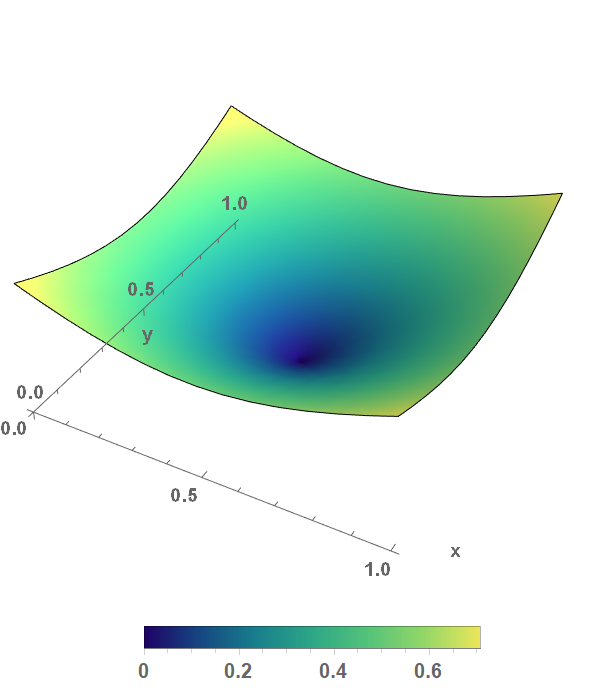}
  \caption{Initial profiles (\ref{quadinit1}) (left) and (\ref{quadinit2}) (right) for example (\ref{rotvelocity}).}
 \label{quadinit}
\end{figure}

We consider two different initial conditions in this example, see Figure \ref{quadinit}. First, the Gaussian is chosen
\begin{equation}
    \label{quadinit1}
    \phi(x,y,0)=e^{-100\left((x-0.5)^2 + (y-0.5)^2\right)} \,,
\end{equation}
and, second, the signed distance function is considered,
\begin{equation}
    \label{quadinit2}
    \phi(x,y,0)=\sqrt{(x-0.5)^2 + (y-0.5)^2} \,.
\end{equation}

In both cases, the initial profile of the solution $\phi$ is deformed up to time $t=0.5$ when the direction of velocity change its sign, so the deformation is reversed afterwards and the initial profile shall be recovered at time $t=1$. 

\begin{figure}[ht]
\centering
 \includegraphics[width=0.49\linewidth]{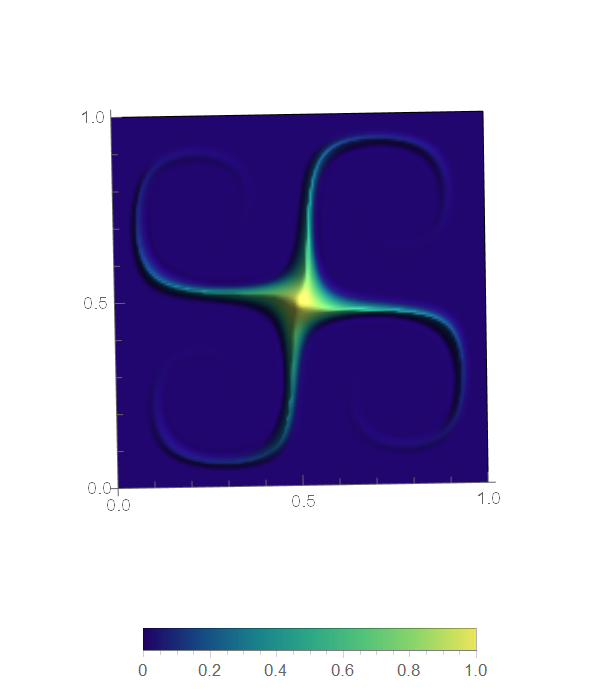}~~
 \includegraphics[width=0.49\linewidth]{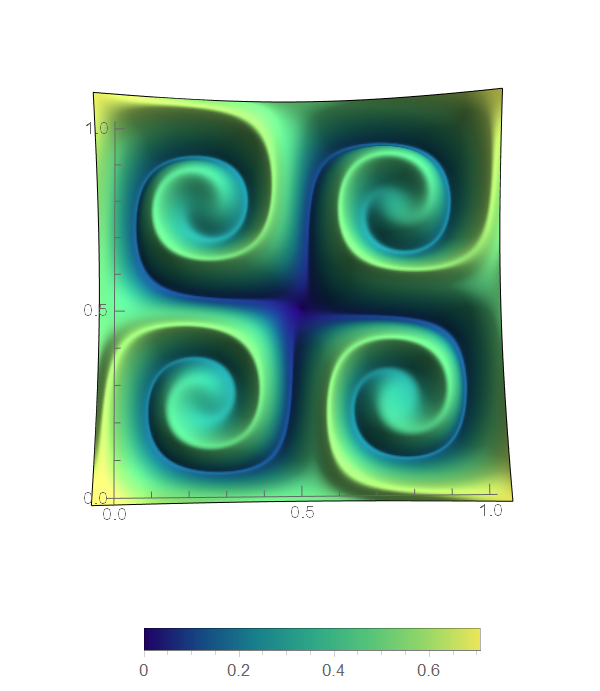}
 \caption{The maximal deformation at $t=0.5$ for (\ref{quadinit1}) (left) and (\ref{quadinit2}) (right) for example (\ref{rotvelocity}).}
 \label{quadhalf}
\end{figure}

We compute the example with maximal Courant numbers being $0.76$ in $y$-direction. The numerical solutions for the finest mesh at the time of maximal deformation are plotted in Figure \ref{quadhalf}. 
%and, moreover, at the final time when the initial profile should be reproduced, see Figure \ref{quadfinal321}.
In Tables \ref{tabquad1} and \ref{tabquad2} we compute the error $E^N$ analogous to (\ref{localerror}) obtained for each numerical solution as the difference between its values at the initial and the final time. 

% \begin{figure}[ht]
% \centering
%  \includegraphics[width=0.49\linewidth]{img/test11_exp_final_321.png}~~
%  \includegraphics[width=0.49\linewidth]{img/test11_sqrt_final_321.png}
%  \caption{Solution at $t=1$ and $I = 320$ for (\ref{quadinit1}) (left) and (\ref{quadinit2}) (right) for example \ref{ex04}}
%  \label{quadfinal321}
% \end{figure}

\begin{table}[ht]
	\begin{center}
	\begin{tabular}{ c c c c c c }
	\hline
 $I$ & $\Delta t$ & $E,\alpha=0.5$ & EOC & $E,\alpha$ in (\ref{kappav}) & EOC \\ 
	\hline
 40 & 100 & 0.01088 & - & 0.00928 & - \\
 80 & 200 & 0.00507 & 1.10 & 0.00415 &  1.16 \\
 160 & 400 & 0.00177 & 1.52 & 0.00138 & 1.59 \\
 320 & 800 & 0.00042 & 2.09 & 0.00030 & 2.18 \\
	\hline 
	\end{tabular}
	\end{center}
	\caption{Numerical error $E^N$ in (\ref{localerror}) for the initial condition (\ref{quadinit1}) in example (\ref{rotvelocity}).}
	\label{tabquad1}
\end{table}

\begin{table}[ht]
	\begin{center}
	\begin{tabular}{ c c c c c c }
	\hline
 $I$ & $\Delta t$ & $E,\alpha=0.5$ & EOC & $E,\alpha$ in (\ref{kappav}) & EOC \\ 
	\hline
 40 & 100 & 0.01692 & - & 0.01355 & - \\
 80 & 200 & 0.00458 & 1.89 & 0.00351 &  1.95 \\
 160 & 400 & 0.00092 & 2.32 & 0.00067 & 2.38 \\
 320 & 800 & 0.00014 & 2.76 & 0.00001 & 2.80 \\
	\hline 
	\end{tabular}
	\end{center}
	\caption{Numerical errors $E^N$ in (\ref{localerror}) with the initial condition (\ref{quadinit2}) in example (\ref{rotvelocity}).}
	\label{tabquad2}
\end{table}

\begin{figure}[ht]
\centering
 \includegraphics[width=0.55\linewidth]{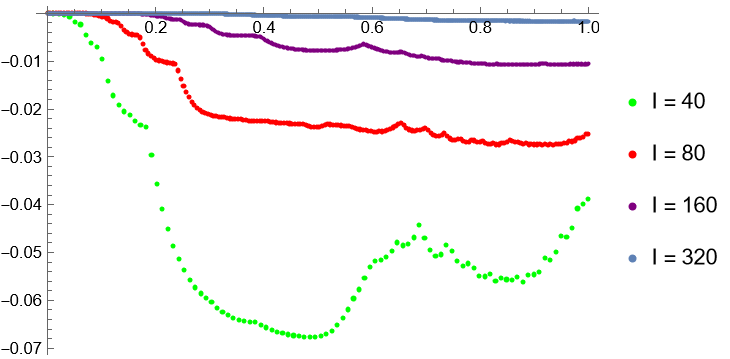}~~
 \includegraphics[width=0.427\linewidth]{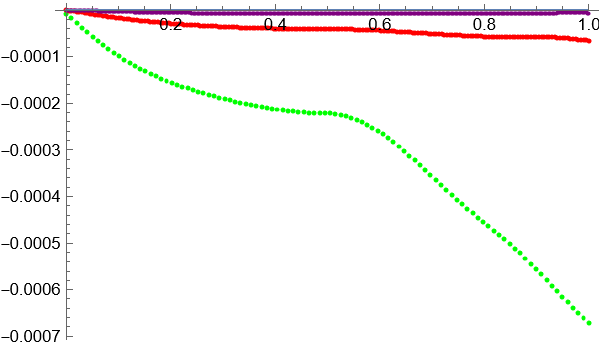}
 
 \caption{The minimum (left) and the difference to the initial mass (right) for numerical solution of example (\ref{rotvelocity}).}
 \label{minquad}
\end{figure}

Next we evaluate the numerical artifact in the form of unphysical negative oscillations in the case of initial condition (\ref{quadinit1}). We plot the minimal value for each numerical solution per each time step at all grid levels in Figure \ref{minquad}. We note that the extremal (rounded) values for each mesh with respect to $n$ are $-0.0677$, $-0.0275$, $-0.0108$, and $-0.00163$ from the coarsest to the finest mesh. 

Finally, the plot of a difference between the initial (conserved) mass and the actual one at each time step for all grid levels can be found in Figure \ref{minquad}. Note that the initial mass (the integral of initial function) is approximated by the value $0.031416$ for (\ref{quadinit1}) and $n=0$ using
\begin{equation}
    \label{mass}
M=M(I,n) = h^2 \sum \limits_{i,j} \phi^n_{i j} \,.
\end{equation}

\section{Conclusions}

We present the novel semi-implicit parametric family of one-dimensional numerical schemes for conservative and non-conservative advection equation. Using the Strang splitting for the advection in several dimensions, one can obtain the numerical solutions of advection equations using in advance known fixed number of alternating substitutions. As the schemes are second order accurate in time and space with unconditional von Neumann stability, they can be considered as a good alternative to standard explicit and implicit schemes for advection dominated problems.

\bibliographystyle{spmpsci}  
\bibliography{main}

\end{document}